\documentclass[11pt]{amsart}

\usepackage{amsmath,amsthm, amscd, amssymb, amsfonts}
\usepackage[all]{xy}
\usepackage{mathrsfs}
\usepackage{enumitem}

\usepackage[T2A,T1]{fontenc}

\usepackage{t1enc}
\usepackage[latin1]{inputenc}

\usepackage{fontsmpl}
\usepackage{srcltx}

\newcommand{\spadechico}{\overset{\spadesuit}{ }}
\newcommand{\starchico}{\overset{\bigstar}{ }}

\allowdisplaybreaks

\newcommand{\supera}[2]{{\mathbf A}(#1|#2)}
\newcommand{\superb}[2]{{\mathbf B}(#1|#2)}

\newcommand{\superd}[2]{{\mathbf D}(#1|#2)}
\newcommand{\superda}[1]{{\mathbf D}(2,1;#1)}
\newcommand{\superf}{{\mathbf F}(4)}
\newcommand{\superg}{{\mathbf G}(3)}

\newcommand{\luq}[1]{\lu_{#1}}

\newcommand{\toba}{{\mathcal B}}

\newcommand{\cO}{{\mathcal O}}
\newcommand{\Vc}{\mathcal{V}}

\newcommand{\sbeta}{\eta}
\newcommand{\sx}{\mathtt x}
\newcommand{\sy}{\mathtt y}
\newcommand{\sxi}{\breve{\xi}}

\newcommand{\ztu}{\overline{\zeta}}

\newcommand{\Dchaintwo}[3]{\xymatrix@C-4pt{\overset{#1}{\underset{\ }{\circ}}\ar
@{-}[r]^{#2}
& \overset{#3}{\underset{\ }{\circ}}}}
\newcommand{\Dchainthree}[5]{\xymatrix@C-6pt{
\overset{#1}{\underset{\ }{\circ}}\ar  @{-}[r]^{#2}  & \overset{#3}{\underset{\
}{\circ}}\ar  @{-}[r]^{#4}
& \overset{#5}{\underset{\ }{\circ}} }}

\newcommand{\Dchainfour}[7]{\xymatrix@C-6pt{\overset{#1}{\underset{\ }{\circ}}\ar
@{-}[r]^{#2}
& \overset{#3}{\underset{\ }{\circ}}\ar  @{-}[r]^{#4}  & \overset{#5}{\underset{\
}{\circ}} \ar  @{-}[r]^{#6}
& \overset{#7}{\underset{\ }{\circ}}}}

\newcommand{\Dchainfive}[9]{\xymatrix@C-6pt{\overset{#1}{\underset{\ }{\circ}}\ar
@{-}[r]^{#2}  & \overset{#3}{\underset{\ }{\circ}}\ar  @{-}[r]^{#4}  &
\overset{#5}{\underset{\ }{\circ}}
\ar  @{-}[r]^{#6}  & \overset{#7}{\underset{\ }{\circ}}\ar  @{-}[r]^{#8}  &
\overset{#9}{\underset{\ }{\circ}}}}

\newcommand{\Dtriangle}[6]{
\xymatrix@R-12pt{  &    \overset{#2}{\circ} \ar  @{-}[dl]_{#4}\ar  @{-}[dr]^{#5} & \\
\overset{#1}{\circ} \ar  @{-}[rr]^{#6}  &  &\overset{#3}{\circ} }}

\newcommand{\Dthreefork}[8]{
\rule[-9\unitlength]{0pt}{12\unitlength}
\begin{picture}(28,12)(0,9)
\put(2,10){\ifthenelse{\equal{#1}{l}}{\circle*{2}}{\circle{2}}}
\put(3,10){\line(1,0){10}}
\put(14,10){\ifthenelse{\equal{#1}{m}}{\circle*{2}}{\circle{2}}}
\put(15,10){\line(1,1){7}}
\put(15,10){\line(1,-1){7}}
\put(22,18){\ifthenelse{\equal{#1}{t}}{\circle*{2}}{\circle{2}}}
\put(22,2){\ifthenelse{\equal{#1}{b}}{\circle*{2}}{\circle{2}}}
\put(2,12){\makebox[0pt]{\scriptsize #2}}
\put(8,11){\makebox[0pt]{\scriptsize #3}}
\put(14,12){\makebox[0pt]{\scriptsize #4}}
\put(19,16){\makebox[0pt][r]{\scriptsize #5}}
\put(19,4){\makebox[0pt][r]{\scriptsize #6}}
\put(24,17){\makebox[0pt][l]{\scriptsize #7}}
\put(24,2){\makebox[0pt][l]{\scriptsize #8}}
\end{picture}}

\newcommand{\Drightofway}[8]{\xymatrix@R-6pt{  &    \overset{#6}{\circ} \ar
@{-}[d]_{#4}\ar  @{-}[dr]^{#7} & \\
\overset{#1}{\circ} \ar  @{-}[r]^{#2}  &\overset{#3}{\circ} \ar  @{-}[r]^{#5}
&\overset{#8}{\circ} }}

\usepackage{color}

\numberwithin{equation}{section}\theoremstyle{plain}

\newtheorem{theorem}{Theorem}[section]

\newtheorem{lemma}[theorem]{Lemma}

\newtheorem{teo intro}{Theorem}

\newtheorem{proposition}[theorem]{Proposition}
\newtheorem{claim}{Claim}

\theoremstyle{definition}

\newtheorem{example}[theorem]{Example}

\theoremstyle{remark}
\newtheorem{remark}[theorem]{Remark}

\newcommand\zt{\Z^{\theta}}

\newcommand\G{\mathbb{G}}

\newcommand\Jb{\mathbb{J}}
\newcommand\I{\mathbb{I}}
\newcommand\id{\operatorname{id}}

\newcommand\ord{\operatorname{ord}}
\newcommand\Hom{\operatorname{Hom}}

\newcommand{\Fr}{\operatorname{Fr}}

\newcommand{\Aut}{\operatorname{Aut}}

\newcommand{\co}{\operatorname{co}}

\newcommand{\ad}{\operatorname{ad}}

\newcommand{\supp}{\operatorname{supp}}
\newcommand{\pma}{\mathbf{p}}
\newcommand{\qmb}{\mathbf{q}}
\newcommand{\ku}{\Bbbk}
\newcommand{\n}{\mathfrak{n}}
\newcommand{\bq}{\mathfrak{q}}
\newcommand{\bp}{\mathfrak{p}}
\newcommand{\bs}{\mathfrak{s}}
\newcommand{\rgo}{\mathfrak{r}}

\newcommand\ot{\otimes}
\newcommand\ra{\rightarrow}

\newcommand\R{\mathbb{R}}
\newcommand\Z{\mathbb{Z}}
\newcommand\N{\mathbb{N}}

\newcommand\cB{\mathcal{B}}
\newcommand\mP{\mathcal{P}}

\newcommand\cU{\mathcal{U}}
\newcommand\cW{\mathcal{W}}
\newcommand\cX{\mathcal{X}}

\newcommand{\az}{\mathfrak{Z}}
\newcommand{\hb}{\mathbf{h}}
\newcommand{\jb}{\mathbf{j}}
\newcommand{\xb}{\mathbf{x}}
\newcommand{\yb}{\mathbf{y}}
\newcommand{\Ht}{\mathtt H}

\newcommand{\br}{\mathfrak{br}}
\newcommand{\brj}{\mathfrak{brj}}
\newcommand{\bgl}{\mathfrak{wk}}
\newcommand{\el}{\mathfrak{el}}
\newcommand{\g}{\mathfrak{g}}
\newcommand{\bg}{\mathfrak{b}}
\newcommand{\ufo}{\mathfrak{ufo}}
\newcommand{\ug}{\mathfrak{u}}

\newcommand\Ufo{\mathtt{ufo}}

\newcommand{\ett}{\mathtt{e}}

\newcommand{\lu}{\mathcal{L}}

\newcommand{\fO}{\mathfrak O}
\newcommand{\wfO}{\underline{\mathfrak O}}
\newcommand{\wbeta}{\underline{\beta}}
\newcommand{\walpha}{\underline{\alpha}}
\newcommand{\wgamma}{\underline{\gamma}}
\newcommand{\wmu}{\underline{\mu}}
\newcommand{\wnu}{\underline{\nu}}

\newcommand{\dpn}{\widetilde{\mathcal{B}}}

\newcommand\pf{\begin{proof}}
\newcommand\epf{\end{proof}}

\begin{document}


\title[Lie algebras arising from Nichols algebras of diagonal type]{Lie algebras arising from Nichols algebras of diagonal type}
\author[Andruskiewitsch; Angiono; Rossi Bertone]
{Nicol\'as Andruskiewitsch, Iv\'an Angiono, Fiorela Rossi Bertone}

\address{FaMAF-CIEM (CONICET), Universidad Nacional de C\'ordoba,
Medina A\-llen\-de s/n, Ciudad Universitaria, 5000 C\' ordoba, Rep\'
ublica Argentina.} \email{(andrus|angiono)@famaf.unc.edu.ar}

\address{Departamento de Matemática, Universidad Nacional del Sur (UNS), Bahía Blanca, Argentina}
\email{fiorela.rossi@uns.edu.ar}

\thanks{The work was partially supported by CONICET,
FONCyT-ANPCyT, Secyt (UNC)}

\begin{abstract}
Let $\mathcal{B}_{\mathfrak{q}}$ be a finite-dimensional Nichols algebra of diagonal type with braiding matrix $\mathfrak{q}$, 
let $\mathcal{L}_{\mathfrak{q}}$ be the corresponding Lusztig algebra as in \cite{AAR1} and let $\operatorname{Fr}_{\mathfrak{q}}: \mathcal{L}_{\mathfrak{q}} \to U(\mathfrak{n}^{\mathfrak{q}})$ 
be the corresponding quantum Frobenius map as in \cite{AAR2}. 
We prove that the finite-dimensional Lie algebra $\mathfrak{n}^{\mathfrak{q}}$ is either 0 or else the positive part of a semisimple Lie algebra $\mathfrak{g}^\mathfrak{q}$
which is determined for each $\mathfrak{q}$ in the list of \cite{H-classif RS}.
\end{abstract}

\maketitle


\section{Introduction}
\subsection{Quantum groups}
Let $\ku$ be an algebraically closed field  of characteristic 0.
The quantized enveloping algebra  of a simple Lie algebra $\g$ was introduced by Drinfeld and Jimbo extending 
the previous definition for $sl(2)$ by Kulish, Reshetikhin and Sklyanin. 
Lusztig introduced and studied the quantum divided power algebra $\cU_q(\g)$ in \cite{L-fdHa-JAMS,L-qgps-at-roots,L-libro}, where $q$ is a root of 1 of odd order (and not divisible by 3 when $\g$ is $G_2$).  It turns out that there is an exact sequence
of Hopf algebras
\begin{align}\label{eq:Lusztig-exact-sequence}
\xymatrix@C=40pt{
\ug_q(\g) \ar@{^(->}[r] & \cU_q(\g) \ar@{->>}[r]^-{\Fr} & U(\g)}
\end{align}
where $\Fr$ was named the quantum Frobenius map by Lusztig and $\ug_q(\g)$ is a finite-dimensional Hopf
algebra, usually called the Frobenius-Lusztig kernel. The exact sequence \eqref{eq:Lusztig-exact-sequence} without
restrictions on $q$ was described in \cite{Lentner}.

In parallel De Concini, Kac and Procesi introduced and studied the quantum group
$U_{q}(\g)$ at a root of unity $q$ (same restrictions as before) in
\cite{DK,DKP,DP}. It is different from $\cU_q(\g)$.
There is an exact sequence of Hopf algebras
\begin{align}\label{eq:DCKP-exact-sequence}
\xymatrix@C=40pt{
\cO (G^d) \ar@{^(->}[r] & U_q(\g) \ar@{->>}[r] & \ug_q(\g)}
\end{align}
where $\cO (G^d)$, a central Hopf subalgebra of $U_{q}(\g)$, is the algebra of functions
on a solvable Poisson algebraic group $G^d$ whose Lie bialgebra is dual to $\g$, 
viewed as Lie bialgebra with the Sklyanin bracket.

The representation theory of $\cU_q(\g)$ has links with algebraic groups in positive
characteristic and with an important class of fusion categories.
As apparent from \eqref{eq:DCKP-exact-sequence} the representation theory of $U_q(\g)$
has a different flavour and is related to the Poisson geometry of $G^d$. 

\subsection{Nichols algebras}
It was recognized early that a key point in the definition of (any version of) quantum groups
is the understanding of the positive part. For those with generic parameter $q$, 
Lusztig and Rosso (and later Schauenburg) gave abstract characterizations of this positive part.
It was then realized that these characterizations fit into frameworks defined independently by Nichols 
(in 1978) and Woronowicz (in 1989).
This abstract notion received the name of Nichols algebra; it is the cornerstone of the method proposed to classify finite-dimensional pointed Hopf algebras \cite{AS-pointed}. 
A basic question is the classification of the finite-dimensional Nichols
algebras of diagonal type. Precisely, let $\bq=(q_{ij})$ be a matrix with entries in $\ku^{\times}$
with connected Dynkin diagram, see \S \ref{subsec:nichols-diag-type}.
It was asked whether 
\begin{align}\label{eq:findim-intro}
\text{the Nichols algebra $\toba_{\bq}$ has  finite dimension.}
\end{align}
The complete answer was provided in  \cite{H-classif RS} through 
the new concepts of (generalized) root systems and  Weyl groupoids \cite{H-Weyl grp,HY}. 
The resulting list 
was organized in various classes in \cite{AA17}: 
\begin{enumerate}[leftmargin=*,label=\rm{(\roman*)}]
\item Cartan and standard types 
(including $\ug_q^+(\g)$ without restrictions on the order of $q$ and twisted versions), 
\item super type (related to Lie superalgebras in characteristic 0), 
\item modular and super modular types (related to Lie algebras and super algebras in characteristic $>0$), and 
\item unidentified type, UFOs for short.
\end{enumerate}

\subsection{Distinguished pre-Nichols algebras}
For the purpose of the classification, it is important to describe the defining relations
of $\toba_{\bq}$ as in \eqref{eq:findim-intro}. This was achieved in \cite{Ang-jems,Ang-crelle} introducing 
along the way  a braided Hopf algebra $\dpn_{\bq}$, further studied in \cite{A-preNichols}.
By \cite[(11)]{AAR2} using \cite[Theorem 29]{A-preNichols} there is an exact sequence
of braided Hopf algebras
\begin{align}\label{eq:intro-pre-Nichols-exact-sequence}
\xymatrix@C=40pt{
Z_{\bq}^+ \ar@{^(->}[r] & \dpn_{\bq} \ar@{->>}[r]^-{\pi_{\bq}} & \toba_{\bq}.}
\end{align}
Assuming the technical condition \eqref{eq:condition cart roots}, $Z_{\bq}^+$ is central in $\dpn_{\bq}$. 
If $\bq$ is of Cartan type and \eqref{eq:condition cart roots} holds, then (up to a twist) $\dpn_{\bq} \simeq U_q^+(\g)$,
the positive part of the quantum group studied in \cite{DK,DKP,DP}. Furthermore
there is a Hopf algebra $U_\bq$ that generalizes the quantum group $U_q(\g)$ and an exact sequence like \eqref{eq:DCKP-exact-sequence} \cite{A-preNichols}.
More generally, if $\bq$ belongs to the \emph{one-parameter family}, then $U_\bq$ is obtained in \cite[Theorem 5.7]{AAY} by evaluation 
at a root of unity of an integral form of a Hopf algebra over (an appropiate localization of ) $\mathbb{C}(q)$. The set of Hopf superalgebras constructed here contains properly the quantized enveloping algebras of basic Lie superalgebras.

\subsection{The Lusztig algebra $\luq{\bq}$}
This is defined as the graded dual of $\dpn_{\bq^t}$ \cite{AAGTV,AAR1}, that in turn is defined by generators and relations \cite{A-preNichols} extending the definition of the De Concini-Kac-Procesi quantum groups \cite{DKP}. At the present moment there is no intrinsic definition of the distinguished pre-Nichols algebra (nor of the De Concini-Kac-Procesi quantum groups). It was conjectured in \cite{A-preNichols} that the distinguished pre-Nichols is the largest covering of the Nichols algebra with finite GK-dim; this is generically true but there are exceptions \cite{ASa,ACS}.

Assume \eqref{eq:condition cart roots}.
If $\bq$ is of Cartan type and \eqref{eq:condition cart roots} holds, then (up to a twist) $\lu_{\bq} \simeq \cU_q^+(\g)$,
the positive part of the quantum group studied in \cite{L-AiM,L-fdHa-JAMS,L-qgps-at-roots}.
Let $\az_\bq$ be the graded dual of $Z_{\bq}^+$ and $\n^\bq=\mP(\az_{\bq})$.
By   \cite{AAR2} $\az_\bq\simeq U(\n^\bq)$
and there is an exact sequence of braided Hopf algebras
\begin{align}\label{eq:extension-braided-lu-intro}
\xymatrix@C=40pt{
\toba_{\bq} \ar@{^(->}[r] & \lu_{\bq} \ar@{->>}[r]^-{\Fr_{\bq}} & \az_{\bq}.}
\end{align}
If $\g$ is a finite-dimensional semisimple Lie algebra  with a Borel subalgebra $\bg$, 
we set $\g_+ := [\bg, \bg]$. 
When $\bq$ has rank 2, the Lie algebra $\n^\bq$ was computed explicitly by a case-by-case analysis \cite{AAR2};
it turns to be either 0 or else isomorphic to $\g^{\bq}_+$, 
where $\g^{\bq}$ is a semisimple Lie algebra (of rank $\leq 2$).
The goal of the present paper is to determine the Lie algebra $\n^\bq$ in general.

\begin{theorem}\label{th:main}  The Lie algebra $\n^\bq$ is either 0 or isomorphic to $\g^{\bq}_+$, 
where $\g^{\bq}$ is a finite-dimensional semisimple Lie algebra as in Tables \ref{table:algebra-cartan},
\ref{table:algebra-super}, \ref{table:algebra-modular}, \ref{table:algebra-super-modular}, \ref{table:algebra-ufo}.
\end{theorem}
In particular we assign semisimple Lie algebras in characteristic 0 to (contragredient)
semisimple Lie (super) algebras in positive characteristic.

The Tables use the notation of \cite{AA17} and refer  to the corresponding section of \emph{loc. cit.}
Reciprocally, the information in the Tables was presented in \cite{AA17} referring to the present paper for proofs.
The Theorem was announced at the H-ACT Conference, Tsukuba,  September 2016 and at the Second Mathematical Congress of the Americas, Montr\'eal, July 2017.

The paper is organized as follows. In Section \ref{sec:preliminaries} we collect definitions and results needed 
for the proof. 
In Section \ref{sec:nilpotent-Lie-algebras} we introduce a set $\wfO^{\bq}$ inside the root lattice
of the Nichols algebra $\toba_{\bq}$ and 
prove that it is a root system as in \cite{Bourbaki}. Let $\g^{\bq}$ be the corresponding semisimple Lie algebra.
The bulk of the proof is Theorem \ref{th:determination-nq} where we show that 
$\n^\bq \simeq \g^{\bq}_+$. In Section \ref{sec:determination} we explain how to compute the Tables.
The Appendix contains a description and a discussion of the algebras appearing in \eqref{eq:Lusztig-exact-sequence},
\eqref{eq:DCKP-exact-sequence}, \eqref{eq:intro-pre-Nichols-exact-sequence} and \eqref{eq:extension-braided-lu-intro}
when $\mathfrak g$ is of type $B_2$.

\smallbreak
The main Theorem is crucial for the development of a Poisson order of the 
Hopf algebra $U_{\bq}$ \cite{A-preNichols}, see \cite{AAY}; it would be useful to study the representation theories of $U_{\bq}$
and of a suitable Drinfeld double of $\az_{\bq}$.

Also, our result contains a uniform description of the positive part of the algebra of \emph{divided powers}, including properly all basic Lie superalgebras. A generalization of \cite[Chapters 34, 35]{L-libro} for super groups is made in \cite{CHW1,CHW2,CSW}, where the authors introduce a parity related with usual root systems; anyway this context does not include new examples of Lie superalgebras except $\mathfrak{osp}(1,2n)$.

\begin{table}
\caption{Cartan and standard types}\label{table:algebra-cartan}
\begin{tabular}{c|c|c}
type of $\bq$ & $\g^\bq$ &  \cite{AA17} 
\\
\hline \hline

 $A_\theta$  & $A_\theta$ & \S 4.1 \\
\hline

 $B_\theta$ & $\begin{matrix}
 B_{\theta},\ N \text{ odd} \\ C_{\theta} ,\ N \text{ even} \end{matrix}$ & \S 4.2 \\
\hline

 $C_\theta$  & $\begin{matrix}
 C_{\theta},\ N \text{ odd} \\ B_{\theta} ,\ N \text{ even} \end{matrix}$ & \S 4.3 \\
\hline

 $D_\theta$  & $D_\theta$ & \S 4.4\\
\hline

 $E_\theta$  & $E_\theta$ & \S 4.5 \\
\hline

 $F_4$   & $F_4$ & \S 4.6 \\
\hline

 $G_2$  & $G_2$ & \S 4.7 \\
\hline

\text{Standard}  $B_\theta$  & $D_k\times  D_{\theta-k}$ & \S 6.1 \\
\hline

\text{Standard}  $G_2$  & $A_1\times  A_1$ & \S 6.2 \\
\hline

\end{tabular}
\end{table}

\bigbreak
\begin{table}
\caption{Super type}\label{table:algebra-super}
\begin{tabular}{c|c|c }
 type of $\bq$ & $\g^\bq$ &  \cite{AA17} 
\\
\hline \hline
  $A(k-1,\theta-k)$ &  $A_{k-1}\times  A_{\theta-k}$ &  \S 5.1 
\\
\hline 
 $B(k,\theta-k)$ &     $\begin{matrix}
C_{k}\times  B_{\theta-k},\ N \text{ odd} \\ C_{k}\times  C_{\theta-k} ,\ N \text{ even} \end{matrix}$ &  \S 5.2
\\
\hline 
 $D(k,\theta-k)$ &   $\begin{matrix}
D_{k}\times  C_{\theta-k},\ N \text{ odd} \\ D_{k}\times B_{\theta-k} ,\ N \text{ even} \end{matrix}$  & \S 5.3
\\ \hline
 $D(2,1,\alpha)$ &    $A_1 \times  A_1\times  A_1$   &  \S 5.4
\\ \hline
 $F(4)$  & 
 $\begin{matrix}
B_3\times  A_1,\ N \text{ odd} \\ C_3\times  A_1 ,\ N \text{ even} \end{matrix}$ &  \S 5.5
\\ \hline 
 $G(3)$  &  $G_2\times  A_1$ &  \S 5.6
\\
\hline 
\end{tabular}
\end{table}

\bigbreak
\begin{table}[ht]
\caption{Modular type}\label{table:algebra-modular}
\begin{tabular}{c|c|c}
type of $\bq$ & $\g^\bq$ & \cite{AA17} \\  
\hline \hline
$\bgl(4,\alpha)$ &   $A_2\times  A_2$ &\S  7.1 
\\ \hline
$\br(2,a)$ &   $A_1\times  A_1$ & \S 7.2
\\ \hline
$\br(3)$ & $B_3$ & \S 7.3 
\\
\hline   
\end{tabular}
\end{table}

\bigbreak
\begin{table}[ht]
\caption{Super modular type 
}\label{table:algebra-super-modular}
\begin{tabular}{c|c|c||c|c|c}
type of $\bq$ & $\g^\bq$ &  \cite{AA17} &  type of $\bq$ &   $\g^\bq$ &  \cite{AA17}
\\
\hline \hline

$\brj(2;3)$ &  $A_1\times  A_1$ & \S 8.1 &  $\g(1,6)$ &  $C_3$ or $A_3$ & \S 8.2
\\
\hline
 $\g(2,3)$ &  $A_2\times  A_1$ & \S 8.3 &
$\g(3,3)$ &  $B_3$ & \S 8.4
\\ \hline
 $\g(4,3)$ &  $C_3\times  A_1$ & \S 8.5
&
$\g(3,6)$ &  $C_4$ & \S 8.6
\\
\hline 
$\g(2,6)$ &  $A_5$ & \S 8.7
& $\el(5,3)$ &  $B_4\times  A_1$ & \S 8.8 
\\
\hline
$\g(8,3)$ &  $F_4\times  A_1$ & \S 8.9 & $\g(4,6)$ &  $D_6$ & \S 8.10 
\\
\hline  $\g(6,6)$ &   $B_6$  & \S 8.11 &$\g(8,6)$ &  $E_7$ & \S 8.12 
\\
\hline $\brj(2,5)$ & $B_2$  & \S 9.1 & $\el(5,5)$ &  $B_5$ & \S 9.2
\\
\hline
\end{tabular}
\end{table}

\bigbreak
\begin{table}[ht]
\caption{Unidentified type 
}\label{table:algebra-ufo}
\begin{tabular}{c|c|c||c|c|c}
type of $\bq$ &   $\g^\bq$ &  \cite{AA17} &  type of $\bq$ &   $\g^\bq$ &  \cite{AA17}
\\
\hline  \hline

$\ufo(1)$ &   $A_5$ & \S 10.1 &  $\ufo(2)$ &   $E_6$ &   \S 10.2
\\
\hline
$\ufo(3)$ & $A_1\times  A_1$ & \S 10.3 & $\ufo(4)$ &  $A_1$ & \S 10.4
\\
\hline 
 $\ufo(5)$  & $A_4$ & \S 10.5 & $\ufo(6)$ & $G_2\times  G_2$ & \S 10.6 
\\
\hline 
$\ufo(7)$ & $0$ & \S 10.7 &  $\ufo(8)$ & $A_1$ & \S 10.8  
\\
\hline 
$\ufo(9)$ & $A_1\times  A_1$ & \S 10.9 &  $\ufo(10)$  & $A_1\times  A_1$ & \S 10.10 
\\
\hline 
$\ufo(11)$ &  $A_1\times  A_1$ & \S 10.11 &  $\ufo(12)$ & $G_2$ & \S 10.12 \\
\hline 
\end{tabular}
\end{table}

\subsection*{Acknowledgements} 
We thank the referees for the careful reading of this paper. We also thank one of them for pointing out Remark \ref{rem:Cartan-root}, which simplifies the proofs of Lemma \ref{lem:cartan-roots-invariant} and Proposition \ref{prop:subalgebra-dist-pre-Nichols} \ref{item:subalgebra-dist-pre-Nichols-i}.

\section{Preliminaries}\label{sec:preliminaries}

\subsection*{Notation}
We fix $\theta \in \N$ and set $\I = \I_{\theta} := \{1, 2, ..., \theta\}$.
Let $(\alpha_j)_{j\in \I}$ be the canonical basis of $\zt$. For each $\beta =\sum_{i\in\I} a_i\alpha_i$, $a_i\in\Z$, the \emph{support} of $\beta$ is 
$\supp \beta= \{i\in\I | a_i\neq 0 \}$.
If $N \in \N$ and $v\in \ku^{\times}$, then $(N)_v := \sum_{j=0}^{N-1}v^{j}$.

For each $q\in\Bbbk^{\times}$, let $\G(q):=\{q^n|n\in\Z\}$.
If $N\in\N$, then $\G_N$ is the group of roots of unity of order $N$ in $\Bbbk$, and $\G_N'$ the subset of primitive roots of order $N$.
Thus $\G(q)=\G_N$ for all $q\in\G_N'$.

\subsection{Nichols algebras of diagonal type}\label{subsec:nichols-diag-type}

Let $\bq=(q_{ij})_{i,j\in\I}$ be a matrix with entries in $\ku^{\times}$.
The Dynkin diagram $\mathfrak{D}$ associated to $\bq$ has set of vertices $\I$, with $i \in \I$ labelled by
$q_{ii}$. Let $i \neq j$ in $\I$ and $\widetilde{q}_{ij} := q_{ij}q_{ji}$.
There is an edge between $i$ and $j$ only when $\widetilde{q}_{ij} \neq 1$, in 
which case this scalar labels the edge: 
$\xymatrix{\overset{q_i}  {\circ} \ar  @{-}[r]^{\widetilde{q}_{ij}}  & \overset{q_j}  {\circ}}$. 
Without loss of generality we shall assume that $\mathfrak{D}$ is connected.
Let $(V,c)$ be a braided vector space of diagonal type associated to $\bq$. Thus, $V$ has a basis
$(x_{i})_{i\in\I}$ such that $c \in GL(V \ot V)$ is given by
\begin{align*}
c(x_i \ot x_j) &= q_{ij} x_j \ot x_i, & i,j&\in\I.
\end{align*}
Let $\cB_{\bq} := \cB(V)$ be the corresponding Nichols algebra (of diagonal type); see \cite{A-leyva,AA17}.
We assume \eqref{eq:findim-intro} i.e. $ \dim \cB_\bq < \infty$. Later we also need the hypothesis 
 \eqref{eq:condition cart roots} that requires more notation. 
 Notice that \eqref{eq:findim-intro} depends only on $\mathfrak{D}$,
 thus it comprises multiparametric versions, while \eqref{eq:condition cart roots} depends fully on $\bq$.

The matrix $\bq$  defines a $\Z$-bilinear form $\qmb:\zt\times\zt\to\ku^\times$ by 
$\qmb(\alpha_j,\alpha_k)=q_{jk}$ for all $j,k\in\I$. We set
\begin{align*}
q_{\alpha\beta} &= \qmb(\alpha,\beta),& \alpha,\beta  &\in \zt.
\end{align*}

The algebra $\cB_{\bq}$ is $\zt$-graded by $\deg x_i := \alpha_i$, $i\in \I$.
By \cite{Kh}, $\cB_\bq$ has a PBW-basis with homogeneous generators.
The \emph{positive} roots of $\cB_{\bq}$ are by definition the elements
of the set $\Delta_+^{\bq}$ of $\zt$-degrees of these generators 
\cite{H-Weyl grp}; it does not depend on the choice of the PBW-basis \cite{H-Weyl grp}. 
The roots of $\cB_{\bq}$  are by definition the elements of
$$\Delta^{\bq} = \Delta_+^{\bq} \cup -\Delta_+^{\bq}.$$

We refer to \cite{AA17} for information on finite-dimensional Nichols algebras of diagonal type.

\subsection{The root system and the Weyl groupoid}
By \eqref{eq:findim-intro},
the following matrix $(c_{ij}^{\bq})_{i,j\in \I}\in\Z^{\theta\times\theta}$ is well-defined \cite{R}:
it is given by  $c_{ii}^{\bq} = 2$ and
\begin{align}\label{eq:defcij}
c_{ij}^{\bq}&:= -\min \left\{ n \in \N_0: (n+1)_{q_{ii}}
(1-q_{ii}^n q_{ij}q_{ji} )=0 \right\},  & i & \neq j.
\end{align}
The matrix $(c_{ij}^{\bq})_{i,j\in \I}$ gives rise to $s_i^{\bq}\in GL(\Z^\theta)$ by 
\begin{align*}
s_i^{\bq}(\alpha_j) &= \alpha_j-c_{ij}^{\bq}\alpha_i,&  j &\in \I,
\end{align*}
$i\in \I$. Clearly $s_i^{\bq}$ is a reflection; it allows to define the matrix
$\rho_i(\bq)$  by
\begin{align}\label{eq:defn-rho-i-q}
\rho_i(\bq)_{jk} &= \qmb(s_i^{\bq}(\alpha_j),s_i^{\bq}(\alpha_k)),& j, k &\in \I.
\end{align}
Let  $\rho_i(V)$ be the braided vector space of diagonal type with matrix $\rho_i(\bq)$.  
Then $\cB_{\bq} \simeq \cB_{\rho_i(\bq)}$ as graded vector spaces \cite{H-Weyl grp}.
Let  
\begin{align}\label{eq:base-groupoid}
\cX := \{\rho_{j_1} \dots \rho_{j_N}(\bq): j_1, \dots, j_N \in \I, N \in \N \}.
\end{align}
We collect the roots $\Delta^{\bp}$ of the matrices
$\bp \in \cX$ as a fibration $\Delta \to  \cX$, with $\Delta^{\bp}$ being the fibre of $\bp$
and call this the root system of $\toba_\bq$.

Let $\cX \times GL(\Z^\theta) \times \cX$ be the groupoid over $\cX$
with source $s$, target $t$ and product given, for $\bp, \rgo, \bs \in \cX$, $g,h \in GL(\Z^\theta)$, by
\begin{align*}
s(\bp, g, \rgo) &= \rgo, &t(\bp, g, \rgo) &= \bp, &
(\bp, g, \rgo)(\bs, h, \bp) &=   (\bs, hg, \rgo).
\end{align*}

The Weyl groupoid $\cW_{\bq}$ of $\cB_\bq$ is the subgroupoid of $\cX \times GL(\Z^\theta) \times \cX$
generated by all 
\begin{align*}
\sigma_{i}^{\bp} &= (\rho_i(\bp), s_i^{\bp}, \bp), &
i &\in \I, \ \bp \in \cX.
\end{align*}
Notice that $\cW_{\bq}$ acts on the root system of $\bq$.  
As shown in \cite{H-Weyl grp},
$\cW_{\bq}$ is finite (conversely, existence and finiteness of $\cW_{\bq}$, together with finitude of the height,  implies $\dim \toba_{\bq} < \infty$). 

\medbreak
Let $w_0 = w_0^{\bq} \in \cW_{\bq}$ be the unique element of maximal length ending in $\bq$ \cite{HY}.
Fix  a reduced expression
\begin{align}\label{eq:reduced-expression}
w_0 &= \sigma_{i_1}^{\bq} \sigma_{i_2}\cdots \sigma_{i_M}.
\end{align}
Then the set $\Delta_+^{\bq}$ can be enumerated as follows, see \cite[Proposition 2.12]{CH1}:
\begin{align} \label{eq:betak}
\Delta_+^{\bq} &= \{\beta_j := s_{i_1}^{\bq}\cdots s_{i_{j-1}}(\alpha_{i_j})\vert \  j\in \I_M\}.
\end{align}
In particular $M$ is the cardinal of the set of positive roots.
We need a generalization of a well-known fact on Weyl groups.

\begin{lemma}\label{lema:reduced-expresion}
Let $i\in\I$, $\bq\in\cX$ and $\bp := \rho_i(\bq)$. 
\begin{enumerate}[leftmargin=*,label=\rm{(\alph*)}]
\item\label{item:reduced-expresion-i} There exists a reduced expression of $w_0 = w_0^{\bq}$ as in \eqref{eq:reduced-expression} such that $i_1=i$.

\item\label{item:reduced-expresion-p} Given the reduced expression with $i_1=i$ as above, there exists $j\in \I$ such that 
$\sigma_{i_2}^{\bp}\dots \sigma_{i_M}\sigma_j$  is a reduced expression of  $w_0^{\bp}$.
\end{enumerate}
\end{lemma}

\pf \ref{item:reduced-expresion-i}:
We prove by induction on $k\le M$ that there exists a family $(i_j)_{j\in\I_k}$ of elements in $\I$ such that $i_1=i$ and $\ell(\sigma_{i_1}^{\bq} \dots \sigma_{i_k}) =k$; the Lemma corresponds to the case $k=M$. The case $k=1$ is evident. Next we assume that $k<M$ and fix
a family $(i_j)_{j\in\I_k}$ such that $i_1=i$  and $\ell(\sigma_{i_1}^{\bq} \dots \sigma_{i_k}) =k$. Set $\bp=\rho_{i_k}\dots \rho_{i_1}(\bq)$. By \cite[Lemma 8 (iii)]{HY},
$k=\vert \Delta_-^{\bq} \cap s_{i_1}^{\bq}\dots s_{i_k}(\Delta_+^{\bp}) \vert$. As $k<M=\vert \Delta_+^{\bp}\vert$, there exists $\beta\in\Delta_+^{\bp}$ such that 
$s_{i_1}^{\bq}\dots s_{i_k}(\beta)\in\Delta_+^{\bq}$. Hence there exists $h\in\supp \beta$ such that $s_{i_1}^{\bq}\dots s_{i_k}(\alpha_h)\in\Delta_+^{\bq}$. We define $i_{k+1} = h$; 
by \cite[Corollary 3]{HY}, $\ell(\sigma_{i_1}^{\bq} \dots \sigma_{i_k}\sigma_{i_{k+1}}) =k+1$.

\ref{item:reduced-expresion-p}:
For, $\sigma_{i_2}^{\bp}\dots \sigma_{i_M}$ is a reduced expression of length $M-1$. 
By \cite[Lemma 8 (iii)]{HY} there exists exactly one positive root 
$\beta\in\Delta_+^{\rho_{i_{M-1}} \dots \rho_{i_2}(\bp)}$ such that 
$s_{i_2}^{\bp}\dots s_{i_M}(\beta)$ is positive. 
Hence $s_{i_2}^{\bp}\dots s_{i_M}(\alpha_j)$ is positive for at least 
one $j \in\supp \beta$, thus $\beta=\alpha_j$. 
By \cite[Corollary 3]{HY}, $\sigma_{i_2}^{\bp}\dots \sigma_{i_M}\sigma_j$ 
is a reduced expression of length $M$, which is $w_0^{\bp}$ by uniqueness.
\epf

The next result about root systems will be useful throughout the article. 

\begin{theorem}{\cite[Theorem 2.4]{CH2}}\label{thm:cuntz-heck}
Let  $\gamma_i\in\Delta_+^{\bq}$, $i\in \I_k$, be linearly independent roots. Then there exist $\bp\in\cX$,  $w\in\Hom (\bq,\bp)$ and $\sigma\in\mathbb{S}_{\I}$ such that the support of $w(\gamma_i)\in\Delta_+^{\bp}$ is contained in $\{\sigma(1),\dots,\sigma(i)\}$ for each $i\in \I_k$. \qed
\end{theorem}

In other words, this result allows to reduce computations on a set of $k$ linearly independent roots to computations on a root system of rank $k$, obtained as a subsystem of a different object of the Weyl groupoid. Clearly 
\begin{align*}
w(\gamma_1)&=\alpha_{\sigma(1)} & \text{ since }\  w(\gamma_1)&\in \Z\alpha_{\sigma(1)} \cap \Delta_+^{\bp} =\{\alpha_{\sigma(1)} \}.
\end{align*}

\medspace

Let $\ug_\bq$ be the Drinfeld double of the bosonization $\cB_{\bq}\# \ku\zt$. The reflections $s_i^\bq$ lift to 
algebra isomorphisms $T^\bq_i : \ug_{\rho_i(\bq)} \ra \ug_\bq$, $i\in\I$,
called Lusztig isomorphisms,
generalizing the isomorphisms of quantized enveloping algebras in \cite{L-libro}; cf. \cite{H-Weyl grp,H-L iso}. 
As for quantum groups, we can then define the root vectors associated to \eqref{eq:reduced-expression} by
\begin{align}\label{eq:root-vectors}
x_{\beta_j}&=T_{i_1}^{\bq}\cdots T_{i_{j-1}}(x_{\alpha_{i_j}}) \in \cB_{\bq}, &\beta_j &\in \Delta_+^{\bq},
\ j\in \I_M.
\end{align}

\subsection{Distinguished-pre-Nichols  algebras}
\subsubsection{Cartan roots}
The notion of Cartan roots  from \cite{A-preNichols} is crucial for the definitions of 
the distinguished pre-Nichols algebra $\dpn_{\bq}$ and its graded dual $\luq{\bq}$.
We say that $i\in\I$ is a \emph{Cartan vertex} of $\bq$  if
\begin{align}\label{eq:cartan-vertex}
q_{ij}q_{ji} &= q_{ii}^{c_{ij}^{\bq}}, & \text{for all } j \neq i.
\end{align}
Then the set of \emph{Cartan roots} of $\bq$ is $\fO^{\bq} = \fO^{\bq}_+ \coprod \fO^{\bq}_-$, where
$\fO^{\bq}_- = -\fO^{\bq}_+$ and
\begin{align*}
\fO^{\bq}_+ &= \big\{s_{i_1}^{\bq} s_{i_2} \dots s_{i_k}(\alpha_i) \in \Delta_+^{\bq}:
i\in \I  \text{ is a Cartan vertex of } \rho_{i_k} \dots \rho_{i_2}\rho_{i_1}(\bq) \big\}.
\end{align*}
Observe that this agrees with the notation in \cite{A-preNichols} but differs from  \cite{AAR1,AAR2}. Set 
\begin{align*}
N^{\bq}_\beta &= \ord q_{\beta\beta} \in\N& &\text{and}& \widetilde N^{\bq}_\beta &= \begin{cases} N^{\bq}_{\beta} &\mbox{ if }\beta\notin\fO^{\bq},
\\ \infty  &\mbox{ if }\beta\in\fO^{\bq}, \end{cases} & \beta &\in \Delta^{\bq}.
\end{align*}

\begin{remark}\label{rem:Cartan-root}
Let $\alpha\in\Delta^{\bq}$. Then $\alpha\in\fO^{\bq}$ if and only if $q_{\alpha\beta}q_{\beta\alpha}\in\G(q_{\alpha\alpha})$ for all $\beta\in\Z^{\theta}$.
\end{remark}
\pf
The case $\alpha=\alpha_j$ is clear by definition of Cartan vertex. By \eqref{eq:defn-rho-i-q}, the image under $s_i^{\bq}$ of the set of $\alpha\in\Z^{\theta}$ such that $q_{\alpha\beta}q_{\beta\alpha}\in\G(q_{\alpha\alpha})$ for all $\beta\in\Z^{\theta}$ is the corresponding set for $\rho_i(\bq)$, so the claim follows.
\epf

\subsubsection{Distinguished-pre-Nichols  algebras}
As in \cite[Definition 1]{A-preNichols} we consider the ideal $\mathcal I_{\bq}$ of $T(V)$ generated by all the defining relations of $\cB_{\bq}$ in \cite[Theorem 3.1]{Ang-crelle}, 
but excluding the power root vectors $x_\beta^{N^{\bq}_\beta}$, $\beta \in\fO^{\bq}_+$,
and adding some quantum Serre relations, see \emph{loc.cit.} for details.
Then $\mathcal I_{\bq}$ is a Hopf ideal \cite[Proposition 3.3]{Ang-crelle} and the distinguished pre-Nichols algebra $\dpn_{\bq}$ of $V$ is defined as the quotient
\begin{align*}
\dpn_{\bq} &= T(V)/\mathcal I_{\bq}.
\end{align*}
Let $Z_{\bq}^+$ be the subalgebra of $\dpn_{\bq}$ generated by $x_{\beta}^{N^{\bq}_{\beta}}$, $\beta \in \fO^{\bq}_+$. Then $Z_{\bq}^+$ is a braided normal Hopf subalgebra of $\dpn_{\bq}$ \cite[Theorem 29]{A-preNichols}. By \cite[(11)]{AAR2}, we have an exact sequence of braided Hopf algebras
\begin{align}\label{eq:pre-Nichols-exact-sequence}
\xymatrix@C=40pt{
Z_{\bq}^+ \ar@{^(->}[r] & \dpn_{\bq} \ar@{->>}[r]^-{\pi_{\bq}} & \toba_{\bq}.}
\end{align}

\medbreak
Let $U_\bq$ be the Drinfeld double of the bosonization $\dpn_{\bq}\# \ku\zt$, see \cite[\S 3, Definition 1]{A-preNichols}.
Thus $U_\bq$ has a triangular decomposition $U_\bq \simeq U_\bq^{+} \otimes U_\bq^{0} \otimes U_\bq^{-}$, where 
$U_\bq^{+} \simeq  \dpn_{\bq}$, $U_\bq^{0} \simeq \ku \Gamma$, $\Gamma \simeq \zt \times \zt$, and $U_\bq^{-}\simeq  \dpn_{\bq^t}$.
We also need the notation $U_\bq^{\leq 0} := U_\bq^{-} \cdot U_\bq^{0} = U_\bq^{0} \cdot U_\bq^{-}$.
There exist
algebra isomorphisms $T^\bq_i : U_{\rho_i(\bq)} \to U_\bq$ 
that descend to the Lusztig isomorphisms $T^\bq_i :\ug_{\rho_i(\bq)} \to \ug_\bq$, $i\in\I$, as above. 

\subsubsection{PBW-basis}
Recall the root vectors $x_{\beta_j}$ in \eqref{eq:root-vectors} associated to the decomposition \eqref{eq:reduced-expression}.
For brevity, $N^{\bq}_j := N^{\bq}_{\beta_j}$, $\widetilde N^{\bq}_j := \widetilde N^{\bq}_{\beta_j}$.
We use below the vector notation
\begin{align*}
\xb^{\hb} &= x_{\beta_M}^{h_M}x_{\beta_{M-1}}^{h_{M-1}} \cdots x_{\beta_1}^{h_1},& &\hb = (h_1,\dots,h_M)\in\N_0^{M}.
\end{align*}
Then  $\dpn_{\bq}$ has a PBW-basis $\Upsilon_{\bq}$ as follows, see \cite[Theorem 11]{A-preNichols}:
\begin{align*}
\Upsilon_{\bq} &=\{ \xb^{\hb}\, | \,  \hb\in\Ht \}, \text{ where}& \Ht &= \{\hb\in\N_0^M: \, 0\leq h_k < \widetilde N^{\bq}_k, \text{ for all } k\in \I_M \}.
\end{align*}
Then 
\begin{align*}
\Upsilon_{\bq} \times \Upsilon_{\bq}  &\simeq \{ \xb^{\hb} \otimes \xb^{\hb'}\, | \,  \hb, \hb' \in\Ht \}&  &\text{is a basis of} & \dpn_{\bq} &\otimes \dpn_{\bq}.
\end{align*}

\subsubsection{Action of $\cW_{\bq}$ on the Cartan roots}
We start  by showing that the sets of Cartan roots are interchanged by the action of the Weyl groupoid.

\begin{lemma}\label{lem:cartan-roots-invariant}
If $i\in\I$, then $N^{\bq}_{\beta} = N^{\rho_i(\bq)}_{s_i^{\bq}(\beta)}$ for all $\beta\in\fO_+^{\bq}$ and $s_i^{\bq}(\fO^{\bq})=\fO^{\rho_i(\bq)}$.
\end{lemma}
\pf
Let $i\in\I$, $\bp=\rho_i(\bq)$ and $\beta\in\fO_+^{\bq}$. Let 
$\pma:\zt\times\zt\to\ku^\times$ be the $\Z$-bilinear form such that
$\pma(\alpha_j,\alpha_k)=p_{jk}$ for all $j,k\in\I$. 
By \eqref{eq:defn-rho-i-q},
\begin{align*}
\pma(\alpha,\beta)&= \qmb(s_i^{\bq}(\alpha),s_i^{\bq}(\beta)) & 
\text{for all }&\alpha,\beta\in\Z^{\I}.
\end{align*}
Hence
\begin{align*}
N^{\rho_i(\bq)}_{s_i^{\bq}(\beta)} = 
\ord \pma(s_i^{\bq}(\beta),s_i^{\bq}(\beta))= 
\ord \qmb(\beta,\beta)= N^{\bq}_{\beta}.
\end{align*}

Finally, the equality $s_i^{\bq}(\fO^{\bq})=\fO^{\bp}$ follows from \eqref{eq:defn-rho-i-q} and Remark \ref{rem:Cartan-root}.
\epf

\subsubsection{Parabolic subalgebras of a distinguished pre-Nichols algebra}
\label{subsubsec:subalgebras-distinguished}

Let $\Jb$ be a subset of $\I$. We identify $\Z^{\Jb}$ with the subgroup of $\Z^{\I}$ generated by $\alpha_j$, $j\in\Jb$. Using this identification, each $\Z^{\Jb}$-graded object is $\Z^{\I}$-graded.

\medspace

Let $W$ be the subspace of $V$ spanned by $x_j$, $j\in\Jb$. Let $\rgo:=(q_{ij})_{i,j\in\Jb}$.
Then $W$ is a braided vector subspace of $V$, with braiding matrix $\rgo$.
The inclusion $W\hookrightarrow V$ induces a map of $\Z^{\I}$-graded Hopf algebras $\varPhi:T(W)\hookrightarrow T(V)$, 
which  descends to an 
injective $\Z^{\I}$-graded Hopf algebra map $\varPhi:\toba_{\rgo}\hookrightarrow \toba_{\bq}$, see \cite[Corollary 2.3]{AS-pointed}.

\begin{remark}\label{rem:root-system-restriction}
The root system of $\rgo$ is the restriction of the root system of $\bq$ to $\Jb$ \cite[Definition 4.1]{CH1}. Indeed $c_{ij}^{\rgo}=c_{ij}^{\bq}$ for all $i,j\in\Jb$, and by \eqref{eq:defn-rho-i-q},
\begin{align*}
\rho_i(\rgo)&=\big(\rho_i(\bq)_{j,k})_{j,k\in\Jb} &&
\text{for all} &&i\in\Jb.
\end{align*}
As the root system of $\bq$ is finite, \cite[Proposition 2.12]{CH1} holds for the restriction to $\Jb$ and this root system coincides with the one of $\rgo$.
\end{remark}

\begin{proposition}\label{prop:subalgebra-dist-pre-Nichols}
\begin{enumerate}[leftmargin=*,label=\rm{(\roman*)}]
\item\label{item:subalgebra-dist-pre-Nichols-i} $\fO_+^{\bq}\cap \Z^{\Jb} \subset \fO_+^{\rgo}$ (here we identify $\fO_+^{\rgo}\subset \Z^{\Jb}$ with the corresponding subset of $\Z^{\I}$).

\smallbreak
\item\label{item:subalgebra-dist-pre-Nichols-ii} The inclusion $W\hookrightarrow V$ induces a $\Z^{\I}$-graded Hopf algebra map 
\begin{align*}
\widetilde{\varPhi}:\dpn_{\rgo}\to \dpn_{\bq}.
\end{align*} 
The image of $\widetilde{\varPhi}$ is the subalgebra $\dpn_{\bq,\Jb}$ of $\dpn_{\bq}$ generated by $x_j$, $j\in\Jb$.

\smallbreak
\item\label{item:subalgebra-dist-pre-Nichols-iii}
 $\widetilde{\varPhi}(Z_{\rgo}^+) = Z_{\bq}^+\cap \dpn_{\bq,\Jb}$ and $\ker \widetilde{\varPhi}_{|Z_{\rgo}^+}$ is the ideal generated by
\begin{align*}
&x_{\beta}^{N_{\beta}}, && \beta\in\fO_+^{\rgo}-\big(\fO_+^{\bq}\cap \Z^{\Jb}\big).
\end{align*}
\end{enumerate}
\end{proposition}

\pf
To prove \ref{item:subalgebra-dist-pre-Nichols-i}, let $\beta\in \fO_+^{\bq}\cap \Z^{\Jb}$. 
By Remark \ref{rem:root-system-restriction}, $\Delta_+^{\bq}\cap\Z^{\Jb}=\Delta_+^{\rgo}$. 
Then the claim follows from Remark \ref{rem:Cartan-root}.

\smallbreak

For \ref{item:subalgebra-dist-pre-Nichols-ii}, we note that each defining relation of $\dpn_{\rgo}$ in \cite[Theorem 3.1]{Ang-crelle} is a defining relation of $\dpn_{\bq}$. Hence the natural map $\varPhi:T(W)\hookrightarrow T(V)$ descends to  $\widetilde{\varPhi}:\dpn_{\rgo}\to \dpn_{\bq}$. As $x_j$, $j\in\Jb$, generate $\dpn_{\rgo}$, the image of $\widetilde{\varPhi}$ is also generated by $x_j$, $j\in\Jb$.

\smallbreak

For \ref{item:subalgebra-dist-pre-Nichols-iii}, let $\pi_{\bq}: \dpn_{\bq}\to \toba_{\bq}$, $\pi_{\rgo}: \dpn_{\rgo}\to \toba_{\rgo}$ be the canonical Hopf algebra projections. 
We consider the following diagram:
\begin{align}\label{eq:pre-Nichols-comm-diagram}
\begin{aligned}
\xymatrix@C=40pt{
Z_{\rgo}^+ \ar@{^(->}[r] \ar@{-->}[d] & \dpn_{\rgo} \ar@{->>}[r]^-{\pi_{\rgo}} \ar@{->}[d]^{\widetilde{\varPhi}} & \toba_{\rgo} \ar@{->}[d]^{\varPhi}
\\
Z_{\bq}^+ \ar@{^(->}[r] & \dpn_{\bq} \ar@{->>}[r]^-{\pi_{\bq}} & \toba_{\bq},}
\end{aligned}
\end{align}
whose rows are the exact sequences \eqref{eq:pre-Nichols-exact-sequence}.
Then the right-hand side square of this diagram is commutative. By \cite[Theorem 29]{A-preNichols}, we have 
$Z_{\rgo}^+=\dpn_{\rgo}^{\co \pi_{\rgo}}$ and $Z_{\bq}^{+} = \dpn_{\bq}^{\co \pi_{\bq}}$, and consequently
\begin{align*}
\widetilde{\varPhi}(Z_{\rgo}^+) = \widetilde{\varPhi}(\dpn_{\rgo}^{\co \pi_{\rgo}}) \subset \dpn_{\bq}^{\co \pi_{\bq}} = Z_{\bq}^+.
\end{align*}
Thus we have a downward map on the left of \eqref{eq:pre-Nichols-comm-diagram} still denoted $\widetilde{\varPhi}$. By  \ref{item:subalgebra-dist-pre-Nichols-ii}, 
$\widetilde{\varPhi}(Z_{\rgo}^+) \subseteq Z_{\bq}^+\cap \dpn_{\bq,\Jb}$.

Next we prove that $\widetilde{\varPhi}(Z_{\rgo}^+) \supseteq Z_{\bq}^+\cap \dpn_{\bq,\Jb}$. We start with a description of $Z_{\bq}^+\cap \dpn_{\bq,\Jb}$.
Let $w_0^{\rgo}=\sigma_{i_1}^{\rgo} \sigma_{i_2} \dots \sigma_{i_L}$ be a reduced expression of the element of maximal length of $\rgo$.
By Remark \ref{rem:root-system-restriction}, $\sigma_{i_1}^{\bq} \sigma_{i_2} \dots \sigma_{i_L}$ is a reduced expression in $\cW_{\bq}$. 
By \cite[Lemma 8 (i)]{HY} there exists a reduced expression $w_0=\sigma_{i_1}^{\bq} \dots \sigma_{i_M}$ of $w_0$ 
that is an extension of $\sigma_{i_1}^{\bq} \sigma_{i_2} \dots \sigma_{i_L}$. Then
\begin{align*}
\Delta_+^{\bq}\cap\Z^{\Jb}=\Delta_+^{\rgo}=\{\beta_t = s_{i_1}^{\bq}\cdots s_{i_{t-1}}(\alpha_{i_t})\vert \  t\in \I_L\}.
\end{align*}
Notice that
\begin{align*}
\dpn_{\bq,\Jb} &= \bigoplus_{\gamma\in(\N_0)^{\Jb}} (\dpn_{\bq})_{\gamma} = \bigoplus_{\gamma\in(\N_0)^{\I}: \, \supp \gamma \subseteq \Jb} (\dpn_{\bq})_{\gamma}.
\end{align*}
By \cite[Proposition 21 \& Theorem 23]{A-preNichols} $Z_{\bq}^+$ is $q$-commutative and has a PBW-basis with generators $x_{\beta}^{N^{\bq}_{\beta}}$, $\beta \in \fO^{\bq}_+$. Hence 
$Z_{\bq}^+\cap \dpn_{\bq,\Jb}$ has a PBW-basis with generators $x_{\beta}^{N^{\bq}_{\beta}}$, $\beta \in \fO^{\bq}_+ \cap\Z^{\Jb}$, so $Z_{\bq}^+\cap \dpn_{\bq,\Jb}$ is the subalgebra generated by $x_{\beta}^{N^{\bq}_{\beta}}$, $\beta \in \fO^{\bq}_+ \cap\Z^{\Jb}$, and is $q$-commutative.

We claim that the following diagram is commutative:
\begin{align}\label{eq:pre-Nichols-Lusztig-iso}
\begin{aligned}
\xymatrix@C=50pt{
U_{\rgo} \ar@{<-}[r]^-{T_{i_1}^{\rgo}} \ar@{->}[d]^{\widetilde{\varPhi}} & 
U_{\rho_{i_1}(\rgo)} \ar@{<-}[r]^-{T_{i_2}^{\rho_{i_1}(\rgo)}} \ar@{->}[d]^{\widetilde{\varPhi}} &
U_{\rho_{i_2}\rho_{i_1}(\rgo)} \ar@{->}[d]^{\widetilde{\varPhi}} \ar@{<.}[r] & \dots
\\
U_{\bq} \ar@{<-}[r]^-{T_{i_1}^{\bq}} & 
U_{\rho_{i_1}(\bq)} \ar@{<-}[r]^-{T_{i_2}^{\rho_{i_1}(\bq)}} &
U_{\rho_{i_2}\rho_{i_1}(\bq)} \ar@{<.}[r] & \dots } 
\end{aligned}
\end{align}
Indeed, by Remark \ref{rem:root-system-restriction} we check recursively on $t$ that the algebra maps 
\begin{align*}
&T_{i_t}^{\rho_{i_{t-1}} \dots \rho_{i_1}(\rgo)}: U_{\rho_{i_{t}} \dots \rho_{i_1}(\rgo)} \to U_{\rho_{i_{t-1}} \dots \rho_{i_1}(\rgo)}, & &\text{and}
\\
&T_{i_t}^{\rho_{i_{t-1}} \dots \rho_{i_1}(\bq)}: U_{\rho_{i_{t}} \dots \rho_{i_1}(\bq)} \to U_{\rho_{i_{t-1}} \dots \rho_{i_1}(\bq)}
\end{align*}
have the same definition on a set of generators of $U_{\rho_{i_{t}} \dots \rho_{i_1}(\rgo)}$ as in \cite[Proposition 3.26]{Ang-crelle}. Using \eqref{eq:pre-Nichols-Lusztig-iso}, $\widetilde{\varPhi}(x_{\beta})=x_{\beta}$ for all $\beta\in \Delta_+^{\bq}\cap\Z^{\Jb}=\Delta_+^{\rgo}$, and
\begin{align}\label{eq:varPhi-generators-image}
\widetilde{\varPhi}(x_{\beta}^{N_{\beta}})=
\begin{cases}
x_{\beta}^{N_{\beta}} & \beta \in \fO^{\bq}_+ \cap\Z^{\Jb}, \\
0 & \beta \notin \fO^{\bq}_+\cap\Z^{\Jb}.
\end{cases}
\end{align}
Then $\widetilde{\varPhi}(Z_{\rgo}^+) = Z_{\bq}^+\cap \dpn_{\bq,\Jb}$. Now $\widetilde{\varPhi}_{|Z_{\rgo}^+}: Z_{\rgo}^+ \to Z_{\bq}^+\cap \dpn_{\bq,\Jb}$ is a map between $q$-polynomial rings. From \eqref{eq:varPhi-generators-image} $\ker \widetilde{\varPhi}_{|Z_{\rgo}^+}$ is the ideal of $Z_{\rgo}^+$ generated by $x_{\beta}^{N_{\beta}}$, $\beta\in\fO_+^{\rgo}-\big(\fO_+^{\bq}\cap \Z^{\Jb}\big)$.
\epf

\begin{remark}\label{rem:subalgebra-dist-pre-Nichols}
The Hopf algebra map $\widetilde{\varPhi}:\dpn_{\rgo}\to \dpn_{\bq}$ of Proposition \ref{prop:subalgebra-dist-pre-Nichols} is not necessarily injective. 

For example, let $\bq$ be a braiding of type $\Ufo(3)$ with Dynkin diagram
\begin{align}\label{eq:dynkin-ufo(3)}
\begin{aligned}
\xymatrix@C35pt{\overset{-1}{\underset{\ }{\circ}}
\ar@{-}[r]^{\ztu}  & 
\overset{\zeta}{\underset{\ }{\circ}}
\ar@{-}[r]^{-\ztu}  &
\overset{-\zeta}{\underset{\ }{\circ}}, }
\end{aligned}
\end{align} 
where $\zeta\in\G_3'$ \cite[\S 10.3]{AA17}. Then 
$\dpn_{\bq}$ is generated by $(x_i)_{i\in \I_3}$ with defining relations
\begin{align}\label{eq:rels-ufo(3)-a}
\begin{aligned}
x_{221}&=0; & x_{13}&=0; & x_1^2&=0; & x_{332}&=0; & x_2^3&=0.
\end{aligned}
\end{align}
Here, $\fO_+^{\bq}=\{\alpha_3, \alpha_1+3\alpha_2+\alpha_3, \alpha_1+3\alpha_2+2\alpha_3 \}$, so $Z_{\bq}^+$ is the subalgebra generated by
\begin{align*}
x_{\alpha_3}^6 &=x_3^6, & 
x_{\alpha_1+3\alpha_2+2\alpha_3}^6 &= [x_{(13)},x_{223}]_c^6, & 
x_{\alpha_1+3\alpha_2+\alpha_3}^6 &= [[x_{(13)},x_2]_c,x_2]_c^6.
\end{align*}

Let $\Jb=\{1,2\}$. Then $\rgo$ is a braiding of type $\supera{0}{1}$ with Dynkin diagram $\xymatrix@C35pt{\overset{-1}{\underset{\ }{\circ}} \ar@{-}[r]^{\ztu}  & 
\overset{\zeta}{\underset{\ }{\circ}} }$ and $\dpn_{\rgo}$ is generated by $(x_i)_{i\in \I_2}$ with defining relations
\begin{align}\label{eq:rels-supera}
x_{221}&=0; & x_1^2&=0.
\end{align}
Here, $\fO_+^{\rgo}=\{\alpha_2\}$, so $Z_{\bq}^+$ is the subalgebra generated by $x_2^3$.
The map $\widetilde{\varPhi}:\dpn_{\rgo}\to \dpn_{\bq}$, $x_i\mapsto x_i$, $i\in\I_2$, annihilates $x_2^3$.
\end{remark}

\subsection{Lusztig algebras}\label{subsec:lusztig-algebras}
The Lusztig algebra $\luq{\bq}$ associated to $\cB_\bq$ is the braided Hopf algebra arising as the 
graded dual of $\dpn_{\bq}$ \cite[3.3.4]{AAGTV}. That is, there is
a bilinear form $\langle \, , \, \rangle: \dpn_{\bq} \times \luq{\bq} \rightarrow \ku$ satisfying
\begin{align*}
\langle y, xx'\rangle &= \langle y_{(2)},x\rangle\langle y_{(1)},x' \rangle, &
\langle yy', x\rangle &= \langle y,x_{(2)}\rangle\langle y',x_{(1)} \rangle,
\end{align*}
for all $x,x' \in \dpn_{\bq}$, $y,y' \in \luq{\bq}$.
Given $\hb\in\Ht$,  let $\yb_{\hb} \in \luq{\bq}$ be determined by 
\begin{align}\label{eq:def-y}
\langle \yb_{\hb},\xb^{\jb}\rangle &= \delta_{\hb, \jb},&  \jb &\in \Ht.
\end{align}

If $k\in \I_M$ and $\hb_k=(0,\dots,1,0\dots,0)\in\N_0^M$ ($1$ in the $k$-th place),  
then we denote $y_{\beta_k}^{(n)} := \yb_{n\hb_k}$.
By \cite[Remark 4.10]{AAR1}, the algebra $\luq{\bq}$ is generated by
\begin{align*}
\{ y_{\alpha_i}: i \in \I\}  \cup  \{y_{\beta}^{(N^{\bq}_\beta)}: \beta\in \fO^{\bq}_+
\text{ such that }
 x_{\beta}^{N^{\bq}_{\beta}}\in\mathcal{P}(\dpn_{\bq}) \}.
\end{align*}

\section{Nilpotent Lie algebras arising from Nichols algebras}\label{sec:nilpotent-Lie-algebras}
\subsection{Extensions}

Let $\az_{\bq}$ be the graded dual of $Z_\bq^+$. 
From now on, we assume the condition
\begin{align}\label{eq:condition cart roots}
&q_{\alpha\beta}^{N^{\bq}_\beta}=1, & \forall \alpha,\beta\in\fO^{\bq}.
\end{align}
As said  \eqref{eq:condition cart roots} depends fully on $\bq$ and not just on the Dynkin diagram. 

\begin{remark}
The assumption \eqref{eq:condition cart roots} implies that $Z_\bq^+$ is \emph{central} in $\dpn_{\bq}$ \cite[Proposition 21]{A-preNichols}.
A similar condition appeared later in \cite{Negron} and ensures the centrality of a sublattice of the weight lattice.
\end{remark}
\bigbreak 
We next recall the nilpotent Lie algebra $\n^\bq$ from \cite{AAR2}.

\begin{proposition}\cite[3.2, 3.3]{AAR2} There is an extension of braided Hopf algebras
\begin{align}\label{eq:extension-braided-lu}
\xymatrix{\cB_{\bq} \ar@{^{(}->}[rr]^{\pi^*}  && \luq{\bq}  \ar@{->>}[rr]^{\iota^*}  &&  \az_{\bq}.}
\end{align}

Let $\n^\bq=\mP(\az_{\bq})$. Then
$\az_\bq$ is a usual Hopf algebra isomorphic to $U(\n^\bq)$ and
$\Big\{\xi_{\wbeta} := \iota^*\Big(y_\beta^{(N^{\bq}_\beta)}\Big): \beta \in \fO^{\bq}_+\Big\}$ is a basis of $\n^\bq$.
\end{proposition}

\subsection{The root system}
In our quest to detemine $\n^\bq$ we  introduce 
\begin{align}\label{eq:root-system-distinguished}
\wfO^{\bq} &= \{N^{\bq}_{\beta}\beta: \beta \in \fO^{\bq} \}, &  \wfO^{\bq}_+ &= \wfO^{\bq} \cap \N_{0}^{\theta},&
\wbeta &= N^{\bq}_{\beta} \beta, \ \beta \in \fO^{\bq}.
\end{align}

\begin{remark}\label{rem:nq-graded} Since $\xi_{\wbeta}$ is homogeneous of degree $\wbeta$,
the Lie algebra $\n^{\bq}$ is $\zt$-graded with support $\wfO_+^{\bq}$ and its components have dimension $\leq 1$.
\end{remark}

\begin{remark}\label{rem:indep-reduced-expression}
The elements $\xi_{\wbeta}$ depend on the choice of the reduced expression \eqref{eq:reduced-expression} of $w_0$. 
However this affects the `root vector' $\xi_{\wbeta}$ only by a non-zero scalar. 
Indeed let $w_0=\sigma_{k_1}^{\bq}\dots \sigma_{k_M}$ be another reduced expression, and
\begin{align*}
\sbeta_j & := s_{k_1}^{\bq}\dots s_{k_{j-1}}(\alpha_{k_j}), & j &\in \I_{M}.
\end{align*}
Then $\big(\sbeta_j\big)_{j \in \I_{M}}$ is just a reordering of $\Delta_+^{\bq}$. We define according to this new reduced expression
as in \eqref{eq:root-vectors}, respectively \eqref{eq:def-y}
\begin{align*}
\sx_{\beta} &\in\dpn_{\bq}, & \sy_{\beta} &\in \luq{\bq}, & \beta\in \Delta_+^{\bq}.
\end{align*}
Let $\beta  \in\fO^{\bq}_+$. 
Then $\sxi_{\wbeta}:=\iota^*\Big((\sy_\beta)^{(N^{\bq}_\beta)}\Big)\in \az_\bq$ satisfies
$\sxi_{\wbeta} \in\Bbbk^{\times} \xi_{\wbeta}$, since both are non-zero elements of degree $\wbeta = N^{\bq}_{\beta}\beta$ in $\n^\bq$.
\end{remark}

Our first  goal is to show that $\wfO^{\bq}$ is a root system inside the $\R$-subspace $\Vc$ of $\R^{\I}$ spanned by $\wfO^{\bq}$. 
For this we follow \cite[VI, D\'efinition 1]{Bourbaki}. By 
construction, $\wfO^{\bq}$ is finite, $0\notin \wfO^{\bq}$ and $\Vc$ is spanned by $\wfO^{\bq}$. Next we need to define for each $\wbeta\in\wfO^\bq$ a suitable element $\wbeta^{\vee}\in \Vc^{\ast}$. We start by some auxiliary results.

\begin{lemma}\label{lem:N-cartan-roots-invariant}
If $i\in\I$, then $s_i^{\bq}(\wfO^{\bq})=\wfO^{\rho_i(\bq)}$.
In particular, $s_i^{\bq}(\wfO^{\bq})=\wfO^{\bq}$ if $i$ is a Cartan vertex.
\end{lemma}
\pf
Let $i\in\I$, $\bp=\rho_i(\bq)$. 
Then $s_i^{\bq}(\wfO^{\bq})=\wfO^{\rho_i(\bq)}$ follows from Lemma \ref{lem:cartan-roots-invariant}. 
The last statement follows since $\rho_i(\bq)=\bq$ if $i$ is a Cartan vertex.
\epf

Let $\wbeta = N^{\bq}_{\beta} \beta\in\wfO^\bq$ and fix an expression
$\beta \overset{\diamond}{=} s_{i_1}^{\bq} \dots s_{i_k}(\alpha_i)$, where $i\in \I$ is a Cartan vertex of $\bp=\rho_{i_k} \dots \rho_{i_1}(\bq)$. We define 
$s_{\wbeta}: \R^{\I}\to \R^{\I}$ by
\begin{align*}
s_{\wbeta} := s_{i_1}^{\bq} \dots s_{i_k} s_i s_{i_k} \dots s_{i_1}.
\end{align*}

\begin{lemma}\label{lem:s-wbeta-properties}
\begin{enumerate}[leftmargin=*,label=\rm{(\alph*)}]
\item\label{item:s-wbeta-properties-i} $\sigma_{i_1}^{\bq}\dots \sigma_{i_k} \sigma_i \sigma_{i_k} \dots \sigma_{i_1} \in\Hom(\bq,\bq)$. Hence $s_{\wbeta}$ is a reflection of $\R^{\I}$ such that
\begin{align*}
s_{\wbeta}(\wbeta) &= -\wbeta, & 
s_{\wbeta}(\Delta^{\bq})&=\Delta^{\bq}, & 
s_{\wbeta}(\wfO^{\bq})&=\wfO^{\bq}.
\end{align*}
\item\label{item:s-wbeta-properties-ii} 
$s_{\wbeta}$ does not depend on $\diamond$: If $\beta=s_{j_1}^{\bq} \dots s_{j_\ell}(\alpha_j)$ is another expression  where $j\in \I$ is a Cartan vertex of $\rgo=\rho_{j_\ell} \dots \rho_{j_1}(\bq)$, then 
$$ s_{\wbeta} = s_{j_1}^{\bq} \dots s_{j_\ell} s_j s_{j_\ell} \dots s_{j_1}. $$
\end{enumerate}
\end{lemma}
\pf
For \ref{item:s-wbeta-properties-i}, notice that $\rho_i(\bp)=\bp$ since $i$ is a Cartan vertex of $\bp$. Hence
\begin{align*}
\rho_i \rho_{i_k} \dots \rho_{i_1}(\bq) 
& = \rho_i (\bp) =\bp,
\end{align*}
so $\sigma_{i_1}^{\bq}\dots \sigma_{i_k} \sigma_i \sigma_{i_k} \dots \sigma_{i_1} \in\Hom(\bq,\bq)$.
This fact implies that
\begin{align*}
s_{\wbeta}(\Delta^{\bq}) &= s_{i_1}^{\bq} \dots s_{i_k} s_i s_{i_k} \dots s_{i_1}(\Delta^{\bq})
= s_{i_1}^{\bq} \dots s_{i_k} s_i (\Delta^{\bp})
\\
&= s_{i_1}^{\bq} \dots s_{i_k} (\Delta^{\bp}) = \Delta^{\bq},
\end{align*}
and similarly $s_{\wbeta}(\wfO^{\bq})=\wfO^{\bq}$, see Lemma \ref{lem:N-cartan-roots-invariant}. 
Set $t=s_{i_1}^{\bq} \dots s_{i_k}\in\Aut \R^{\I}$, so $s_{\wbeta}=t s_i^{\bp}t^{-1}$. Now $s_{\wbeta}$ is a reflection since $s_i^{\bp}$ is a reflection, and 
\begin{align*}
s_{\wbeta}(\wbeta) &=N^{\bq}_{\beta} s_{\wbeta}(\beta) = N^{\bq}_{\beta} t s_i t^{-1} (t(\alpha_i))=- N^{\bq}_{\beta} t(\alpha_i)= -N^{\bq}_{\beta} \beta =-\wbeta.
\end{align*}
For \ref{item:s-wbeta-properties-ii}, set $s:= s_{j_1}^{\bq} \dots s_{j_\ell} s_j s_{j_\ell} \dots s_{j_1}: \R^{\I}\to \R^{\I}$. 
By \ref{item:s-wbeta-properties-i}, $s$ is a reflection such that $s(\beta)=-\beta$ and $s(\Delta^{\bq})=\Delta^{\bq}$, hence $s=s_{\wbeta}$ by \cite[IV, Lemme 1]{Bourbaki}.  
\epf

As $s_{\wbeta}$ preserves $\wfO^\bq$ and $s_{\wbeta}(\wbeta) = -\wbeta$, 
the restriction of $s_{\wbeta}$ is a reflection of $\Vc$ (that we denote by the same name). 
Thus we define $\wbeta^{\vee} \in \Vc^{\ast}$ as the unique linear form such that
\begin{align*}
s_{\wbeta}(v) &= v- \wbeta^{\vee}(v) \wbeta & & \text{for all }v\in \Vc.
\end{align*}

The strategy to reach our goal in the proof of the next Theorem is to reduce to  rank 2 and invoke \cite{AAR2}.

\begin{theorem}\label{thm:cartan-root-system}
$\wfO^{\bq}$ is a root system of $\Vc$.
\end{theorem}
\pf By \cite[VI, D\'efinition 1]{Bourbaki} and  Lemma \ref{lem:s-wbeta-properties}, it only remains to show that 
\begin{align*}
\wbeta^{\vee}(\wfO^{\bq}) &\subset \Z & \text{ for all }  \wbeta &\in \wfO^{\bq}.
\end{align*}
Let $\wbeta,\wgamma\in\wfO^{\bq}$, where $\wbeta=N^{\bq}_{\beta}\beta$, 
$\wgamma=N^{\bq}_{\gamma}\gamma$, $\beta\neq \gamma\in\Delta_+^{\bq}$. By Theorem \ref{thm:cuntz-heck}, there exist $\bp\in\cX$, $w\in\Hom(\bq,\bp)$ and $i,j\in\I$ such that
\begin{align*}
w(\beta)&=\alpha_i, & \alpha:=w(\gamma)&\in\Delta^{\bp} \cap (\N_0\alpha_i+\N_0\alpha_j).
\end{align*}
Notice that $N^{\bp}_i:=\ord p_{ii}$ coincides with $N^{\bq}_{\beta}$, and also $N^{\bq}_{\gamma}=N^{\bq}_{\alpha}$, so $\walpha=w(\wgamma)$. 
Let $\rgo$ be the $2\times 2$-submatrix of $\bp$ corresponding to $i,j$. Then 
$$\Delta^{\rgo} = \Delta^{\bp}\cap(\Z\alpha_{i}+\Z\alpha_j).$$ 
By \cite{AAR2}, $\n^{\rgo}$ is the positive part of a (rank 2) semisimple Lie algebra, 
so $\wfO^{\rgo}$ is a root system. As $s_{i}^{\rgo}(\walpha)\in\wfO^{\rgo}$, there exists $b\in\Z$ such that $s_i^{\rgo}(\walpha)=\walpha-b\walpha_i$. As ${s_i^{\bp}}_{|\Z\alpha_{i}+\Z\alpha_j}$ depends only on $c_{ij}^{\bp}=c_{ij}^{\rgo}$, we have that ${s_i^{\bp}}_{|\Z\alpha_{i}+\Z\alpha_j}=s_i^{\rgo}$: in particular, 
$s_i^{\bp}(\walpha)=s_i^{\rgo}(\walpha)=\walpha-b\walpha_i$.

We fix a reduced expression $w^{-1}=\sigma_{i_1}^{\bq}\dots \sigma_{i_k}$. Let $t=s_{i_1}^{\bq}\dots s_{i_k}$. As $\wbeta=w^{-1}(\alpha_i)$, we have that $s_{\wbeta}=ts_{i}^{\bp}t^{-1}$, see Lemma \ref{lem:s-wbeta-properties} \ref{item:s-wbeta-properties-ii}. Hence
\begin{align*}
s_{\wbeta}(\wgamma) &=ts_{i}^{\bp}t^{-1}(\wgamma) = ts_i^{\bp}(\walpha) = t(\walpha-b\walpha_i)=\wgamma-b\wbeta.
\end{align*}
Thus $\wbeta^{\vee}(\wgamma)=b\in\Z$.
\epf

The following is a set of simple roots of $\wfO^{\bq}$:
\begin{align}\label{eq:Oq-simple}
\varPi^{\bq} &= \{\varpi \in \wfO^\bq_+: \varpi \neq \alpha + \beta \mbox{ for all } \alpha,\beta \in \wfO^\bq_+ \}.
\end{align}

We next prove that the root system is the same for all $\bp \in \cX$, see \eqref{eq:base-groupoid}.
(Recall that such $\bp$ is said to be Weyl equivalent to $\bq$).

\begin{lemma}\label{lemma:root-system-weyl-equivalence}
Let $\bp\in\cX$. Then $\wfO^{\bq}$ and $\wfO^{\bp}$ are isomorphic as root systems.
\end{lemma}
\pf
It is enough to consider the case $\bp=\rho_i(\bq)$, $i\in\I$. Let $\Vc'$ be the $\R$-subspace of $\R^{\I}$ spanned by $\wfO^{\bp}$. By Lemma \ref{lem:N-cartan-roots-invariant} $s_i^{\bq}:\Vc\to\Vc'$ is a linear isomorphism such that $s_i^{\bq}(\wfO^{\bq})=\wfO^{\bp}$. It remains to prove that
\begin{align*}
s_{s_i^{\bq}(\wbeta)}&= s_i^{\bq} s_{\wbeta} (s_i^{\bq})^{-1}
=s_i^{\bq} s_{\wbeta} s_i^{\bq} & 
\text{for all }& \wbeta\in \wfO^{\bq}.
\end{align*}
Let $\wbeta\in \wfO^{\bq}$, $\wbeta'=s_i^{\bq}(\wbeta) \in \wfO^{\bp}$. Let $\beta\in\fO^{\bq}$ be such that $\wbeta=N_{\beta}\beta$. If $\beta=\alpha_i$, then $s_{\wbeta} =s_i^{\bq}=s_i^{\bp}=s_{\wbeta'}$.

Assume now $\beta\neq\alpha_i$.
By definition there exist $i_1,\dots,i_k,j\in\I$ such that $\beta'=s_{i_1}^{\bp} \dots s_{i_k}(\alpha_j)$, where $j$ is a Cartan vertex of $\rgo=\rho_{i_k} \dots \rho_{i_1} (\bp)$. Hence $\beta=s_i^{\bq} s_{i_1} \dots s_{i_k}(\alpha_j)$: By Lemma \ref{lem:s-wbeta-properties} \ref{item:s-wbeta-properties-ii},
\begin{align*}
s_{\wbeta} &= s_i^{\bq} s_{i_1} \dots s_{i_k} s_j s_{i_k} \dots s_{i_1} s_i,
\\
s_{\wbeta'} &= s_{i_1}^{\bp} \dots s_{i_k} s_j s_{i_k} \dots s_{i_1} = (s_i^{\bq})^2 s_{i_1}^{\bp} \dots s_{i_k} s_j s_{i_k} \dots s_{i_1} (s_i^{\bq})^2 = s_i^{\bq} s_{\wbeta} s_i^{\bq}.
\end{align*}
Thus $s_i^{\bq}$ is an isomorphism of root systems.
\epf

\subsection{Identification of $\n^{\bq}$} Let $\g^{\bq}$ be the semisimple Lie algebra with root system $\wfO^\bq$, let $\bg^{\bq}$ be a  
Borel subalgebra and $\g^{\bq}_+ = [\bg^{\bq}, \bg^{\bq}]$. Here is our next goal:

\begin{theorem}\label{th:determination-nq} There is an isomorphism of Lie algebras
$\n^\bq \simeq \g^{\bq}_+$.
\end{theorem}

For this, we use the (probably well-known) characterization:

\begin{lemma}\label{lem:graded-Lie-algebra-positive-roots}
Let $\g_+$ be the positive part of a semisimple Lie algebra $\g$ of rank $r$ with root system $\Delta$
and root lattice $Q$.
Let $\n=\oplus_{\beta \in Q} \n_{\beta}$ be a $Q$-graded Lie algebra such that:

\begin{enumerate}[leftmargin=*,label=\rm{(\alph*)}]
\item\label{item:graded-Lie-algebra-positive-roots-i} For each $\beta\in Q$, $\dim \n_{\beta}= \begin{cases}
1 & \text{if }\beta\in\Delta_+, \\ 0 & \text{if }\beta\notin\Delta_+.
\end{cases}$
\item\label{item:graded-Lie-algebra-positive-roots-ii} If $\mu, \nu \in\Delta_+$ are such that $\mu+\nu \in\Delta_+$, then $\n_{\mu+\nu}=[\n_{\mu} , \n_{\nu}]$.
\end{enumerate}
Then $\n \simeq \g_+$ as Lie algebras.
\end{lemma}

\pf
Let $C = (c_{ij})$ be the Cartan matrix of $\g$. Recall that $\g_+$ is presented by generators $e_i$, $i\in\I$, and the Serre relations
$(\ad e_i)^{1-c_{ij}} e_j=0$ for $j\neq i$. 

For each $\beta\in\Delta_+$ we pick $\ett_{\beta}\in\n_{\beta}-0$, hence $\n_{\beta}= \Bbbk \ett_{\beta}$ by \ref{item:graded-Lie-algebra-positive-roots-i}. First we claim that there exists a Lie algebra map $\varphi:\g_+\to \n$ such that
$\varphi(e_i)=\ett_{\alpha_i}$ for all $i\in\I$. In fact, for $i\neq j\in\I$, $(\ad \ett_i)^{1-c_{ij}} \ett_j \in \n_{(1-c_{ij})\alpha_i+\alpha_j}$. By \ref{item:graded-Lie-algebra-positive-roots-i} $\n_{(1-c_{ij})\alpha_i+\alpha_j}=0$, since $(1-c_{ij})\alpha_i+\alpha_j\notin \Delta_+$. Hence 
$(\ad \ett_i)^{1-c_{ij}} \ett_j=0$ for all $j\neq i$, which implies the existence of $\varphi$.

Next we claim that $\varphi$ is surjective. By \ref{item:graded-Lie-algebra-positive-roots-i} it suffices to prove that $\ett_{\beta}\in \varphi(\g_+)$ for all $\beta\in\Delta_+$. Let $\beta\in\Delta_+$: by \cite[VI, Proposition 19]{Bourbaki} we may write $\beta=\alpha_{i_1}+\dots+\alpha_{i_n}$, $i_j\in\I$, such that
each $\gamma_k:= \alpha_{i_1}+\dots+\alpha_{i_k}\in\Delta_+$ for all $k\in\I_n$. By \ref{item:graded-Lie-algebra-positive-roots-ii}, 
$\ett_{\gamma_{k+1}}\in\Bbbk [\ett_{\gamma_k}, \ett_{\alpha_{i_{k+1}}}]$ for all $k<n$. Hence we have that $\ett_k\in \varphi(\g_+)$ by induction on $k$; in particular $\beta=\gamma_n\in \varphi(\g_+)$.

Finally $\varphi$ is an isomorphism since $\dim \n= \vert \Delta_+ \vert = \dim \g_+$ by \ref{item:graded-Lie-algebra-positive-roots-i}.
\epf

We have seen in Remark \ref{rem:nq-graded} that  $\n^{\bq}$ satisfies hypothesis
\ref{item:graded-Lie-algebra-positive-roots-i}  of Lemma \ref{lem:graded-Lie-algebra-positive-roots}; 
we next prove that also \ref{item:graded-Lie-algebra-positive-roots-ii} holds.
For this, we reduce again to rank 2, where is already known \cite{AAR2},
moving around by the action of the Weyl groupoid.

\begin{lemma}\label{lemma:generators-nq}
Let $\wmu$, $\wnu$, $\wgamma \in \wfO^\bq_+$ be such that
$\wmu + \wnu = \wgamma$. Then 
\begin{align}\label{eq:alpha-beta-gamma}
[\xi_{\wmu},\xi_{\wnu}] &\neq 0.
\end{align}
Consequently, the Lie algebra $\n^{\bq}$ is generated by 
$\varPi^{\bq}$.
\end{lemma}

\pf 
To start with, choose $i\in \I$ and set $\bp=\rho_i(\bq)$, 
$\mu' := s_i(\mu)$, 
$\nu' := s_i(\nu)$, $\gamma' :=s_i(\gamma)$.
By Lemma \ref{lem:N-cartan-roots-invariant}, we have
\begin{align*}
\wmu' &= s_i(\wmu), & \wnu' &= s_i(\wnu),&  \wgamma' &= s_i(\wgamma) \in \wfO^{\bp}.
\end{align*}

\begin{claim}\label{claim:uno} If   $\wmu'$, 
$\wnu'$, $\wgamma'\in \wfO^{\bp}_+$ and $[\xi_{\wmu'},\xi_{\wnu'}] \neq 0$ in $\n^{\bp}$, 
then $[\xi_{\wmu},\xi_{\wnu}] \neq 0$ in $\n^{\bq}$.
\end{claim}
\noindent \emph{Proof of the Claim.}
By Lemma \ref{lema:reduced-expresion} there is a reduced expression
$$w_0 \overset{\bigstar}{=} \sigma_{i_1}^{\bq} \dots \sigma_{i_M}$$ 
with $i=i_1$ and $h\in\I$ 
such that $w_0^{\bp} \overset{\spadesuit}{=} \sigma_{i_2}^{\bp} \dots \sigma_{i_M} \sigma_h$ is a reduced expression of $w_0^{\bp}$. 
By Remark \ref{rem:indep-reduced-expression} we may assume that the root vectors 
 $x_{\beta}$, $\beta\in \Delta_+^{\bq}$, for $\bq$ are defined using the reduced expression $\starchico$. 
As in \eqref{eq:betak} we set 
\begin{align*}
\beta_t &= s_{i_1}^{\bq}\dots s_{i_{t-1}}(\alpha_{i_t}),& t &\in \I_M.
\end{align*}
Assume that $\mu<\nu$ for the convex order associated to $\starchico$. 
Hence $\mu<\gamma<\nu$ and there exists $j<k<\ell\in\I_{2,M}$ such that
$\mu= \beta_{j}$, $\gamma = \beta_{k}$, $\nu = \beta_{\ell}$, see \cite[Theorem 3.11]{Ang-jems}.
(Here $j >1$ since $\beta_1 = \alpha_{i_1}= \alpha_{i}$ and then $s_i(\beta_1) = -\alpha_{i} \in \Delta^{\bp}_-$).
We enumerate $\Delta_+^{\bp}$ using 
$\spadechico$ as in \eqref{eq:betak}:
\begin{align*}
\beta'_t &= s_{i_2}^{\bq}\dots s_{i_{t-1}}(\alpha_{i_t}), 
\quad t\in \I_{2,M}, & 
\beta_{M+1}' &= s_{i_2}^{\bq}\dots s_{i_{M}}(\alpha_{i_h})=\alpha_i.
\end{align*}
Therefore, $s_i(\beta'_t)=\beta_t$ for all $t \in \I_{2,M}$, so
\begin{align*}
\mu' &= s_{i_2}^{\bp}\dots s_{i_{j-1}}(\alpha_{i_j}), &
\gamma' &= s_{i_2}^{\bp}\dots s_{i_{k-1}}(\alpha_{i_k}), &
\nu' &= s_{i_2}^{\bp}\dots s_{i_{l-1}}(\alpha_{i_l}).
\end{align*}
By Remark \ref{rem:indep-reduced-expression} again, the root vectors 
$x_{\beta'}$, $\beta'\in \Delta_+^{\bp}$, for $\bp$ are defined 
using the reduced expression $\spadechico$; one should not confuse them with  the $x_{\beta}$'s
that correspond to $\bq$.
Using this, we see by definition that 
\begin{align}\label{eq:Ti-root-vectors}
T_i^{\bq}(x_{\beta'_t})&=x_{\beta_t} \quad \text{for all }t \in\I_{2,M}, 
& T_i^{\bq}(x_{\beta'_{M+1}})& \in U_\bq^{\le 0}.
\end{align}
In particular we have that
\begin{align*}
T_i^{\bq}(x_{\mu'}) &= x_{\mu},& T_i^{\bq}(x_{\nu'}) &= x_{\nu},& T_i^{\bq}(x_{\gamma'}) &= x_{\gamma}.
\end{align*}

As $\wmu + \wnu = \wgamma$, there exists $c\in\ku$ such that $[\xi_{\wmu},\xi_{\wnu}]=c\xi_{\wgamma}$;
our goal is to prove that $c \neq 0$. 
Let $\varDelta$ be the comultiplication of $\dpn_\bq$. Then
\begin{align*}
c = c \, \xi_{\wgamma}(x_{\gamma}^{N_{\gamma}}) =[\xi_{\wmu},\xi_{\wnu}] (x_{\gamma}^{N_{\gamma}})
=(\xi_{\wmu}\ot\xi_{\wnu}-\xi_{\wnu}\ot\xi_{\wmu}) \varDelta(x_{\gamma}^{N_{\gamma}}).
\end{align*}
By \cite[Proposition 28]{A-preNichols}, we have
\begin{align}\label{eq:coproduct-xgamma-left}
\varDelta(x_{\gamma}^{N_{\gamma}}) = \varDelta(x_{\beta_k}^{N_{k}}) & 
\in x_{\beta_k}^{N_{k}} \otimes 1 + \sum   x_{\beta_{k-1}}^{h_{k-1}}  x_{\beta_{k-2}}^{h_{k-2}} \dots  x_{\beta_{1}}^{h_{1}} \otimes \dpn_{\bq}.
\end{align}
Recall that $x_i=x_{\beta_1}$. Then $x_{\nu}^{N^{\bq}_{\nu}} \ot x_{\mu}^{N^{\bq}_{\mu}} = x_{\beta_{\ell}}^{N_{\ell}} \otimes x_{\beta_j}^{N_{j}}$ 
does not appear in this expression of $\varDelta(x_{\gamma}^{N_{\gamma}})$. We conclude by \eqref{eq:def-y} that
\begin{align*}
-c = (\xi_{\wnu}\ot\xi_{\wmu}) \varDelta(x_{\gamma}^{N_{\gamma}})
\end{align*}
is the coefficient of $x_{\mu}^{N^{\bq}_{\mu}}\ot x_{\nu}^{N^{\bq}_{\nu}}$  in the expression of $\varDelta(x_{\gamma}^{N^{\bq}_{\gamma}})$ in the basis $\Upsilon_{\bq} \times \Upsilon_{\bq}$; so we have to prove that this is $\neq 0$. Thus Claim \ref{claim:uno}
is equivalent to:

\begin{claim} \label{claim:dos}
If the coefficient of  $x_{\mu'}^{N^{\bp}_{\mu'}}\ot x_{\nu'}^{N^{\bp}_{\nu'}}$ in the expression of $\varDelta(x_{\gamma'}^{N_{\gamma'}})$ 
in the basis $\Upsilon_{\bp} \times \Upsilon_{\bp}$ is $\neq 0$, then so is the coefficient of $x_{\mu}^{N^{\bq}_{\mu}}\ot x_{\nu}^{N^{\bq}_{\nu}}$ in the expression of $\varDelta(x_{\gamma}^{N^{\bq}_{\gamma}})$ in the basis $\Upsilon_{\bq} \times \Upsilon_{\bq}$.
\end{claim}

To deal with Claim \ref{claim:dos}, we need more facts from \cite{A-preNichols}. As in \cite[Remark 5]{A-preNichols}, where more details could be found, 
let
\begin{itemize}[leftmargin=*]\renewcommand{\labelitemi}{$\circ$}
\item $\dpn_{\bp,i} = \ku[x_i]$, the subalgebra of $\dpn_{\bp}$ generated by $x_i$, 

\smallbreak
\item $\iota_i: \dpn_{\bp,i} \to \dpn_{\bp}$ the inclusion,

\smallbreak
\item $\pi_i:\dpn_{\bp}\to \dpn_{\bp,i}$ the projection annihilating $x_j$ for all $j\neq i$,

\smallbreak
\item  $\triangleleft$ the action of $(\dpn_{\bp,i})^*$ on $\dpn_{\bp}$ (whose precise definition is not needed here), 

\smallbreak
\item and $y_i^{(k)}\in (\dpn_{\bp,i})^*$ such that $y_i^{(k)}(x_i^l)=\delta_{k,l}$. Notice that $y_i^{(k)}$ is nothing but the image of $y_{\alpha_{i}}^{(k)}$ in \S \ref{subsec:lusztig-algebras} under the map ${\iota_i}^t: (\dpn_{\bp})^* \to (\dpn_{\bp,i})^*$.
\end{itemize}

By \cite[Theorem 27]{A-preNichols} we have
\begin{align}\label{eq:coprod-x_gamma}
\begin{aligned}
&\varDelta(x_{\gamma}^{N^{\bq}_\gamma}) =
\varDelta\Big(T_i^{\bq}(x_{\gamma'}^{N^{\bp}_{\gamma'}}) \Big)
\\ &= \sum_{k\geq0} T_i^{\bq} \Big( (x_{\gamma'}^{N^{\bp}_{\gamma'}})_{(1)} \Big) x_i^{k} \otimes T_i^{\bq} \Big(\iota_i\mathcal S^{-1}\pi_i \big((x_{\gamma'}^{N^{\bp}_{\gamma'}})_{(2)}\big) \big(x_{\gamma'}^{N^{\bp}_{\gamma'}}\big)_{(3)} \Big)    \triangleleft y_i^{(k)}.
\end{aligned}
\end{align}

To compute the coefficient of $x_{\mu}^{N^{\bq}_\mu} \otimes x_{\nu}^{N^{\bq}_\nu}$ in \eqref{eq:coprod-x_gamma}, written as a linear combination of elements of the PBW-basis of $\dpn_{\bq}$ in both sides of the tensor product, we write each term of
\begin{align}\label{eq:xgammaprima-N-PBW}
(\id \ot \varDelta) \varDelta(x_{\gamma'}^{N^{\bp}_{\gamma'}})= (x_{\gamma'}^{N^{\bp}_{\gamma'}})_{(1)} \otimes  (x_{\gamma'}^{N^{\bp}_{\gamma'}})_{(2)} \ot   (x_{\gamma'}^{N^{\bp}_{\gamma'}})_{(3)}
\end{align}
as a linear combination of the elements of the PBW-basis of $\dpn_{\bp}$. 
By \cite[Proposition 28]{A-preNichols}, now applied to $x_{\gamma'}^{N_{\gamma'}}=x_{\beta'_k}^{N_{\beta'_k}} \in \dpn_{\bp}$,
\begin{align*}
(\id \ot \varDelta) \varDelta(x_{\gamma'}^{N_{\gamma'}}) 
& \in x_{\beta'_k}^{N_{k}} \otimes 1 \otimes 1 
+ \sum   x_{\beta'_{k-1}}^{h_{k-1}}  x_{\beta'_{k-2}}^{h_{k-2}} \dots  x_{\beta'_{2}}^{h_{2}} \otimes \dpn_{\bp}\otimes \dpn_{\bp}.
\end{align*}
Using \eqref{eq:Ti-root-vectors} again, $x_{\mu}^{N^{\bq}_\mu} \otimes x_{\nu}^{N^{\bq}_\nu}$ can appear in 
\begin{align*}
T_i^{\bq} \left( (x_{\gamma'}^{N^{\bp}_{\gamma'}})_{(1)} \right) x_i^{k} \otimes T_i^{\bq} \left(\iota_i\mathcal S^{-1}\pi_i \Big((x_{\gamma'}^{N^{\bp}_{\gamma'}})_{(2)}\Big) \Big(x_{\gamma'}^{N^{\bp}_{\gamma'}}\Big)_{(3)} \right)    \triangleleft y_i^{(k)}
\end{align*}
only when $k=0$ and for those terms of \eqref{eq:xgammaprima-N-PBW} with $x_{\mu'}^{N^{\bp}_{\mu'}}$ on the left hand side. 

Next we determine the middle part of those terms of \eqref{eq:xgammaprima-N-PBW} that contribute to the coefficient
of $x_{\mu}^{N^{\bq}_\mu} \otimes x_{\nu}^{N^{\bq}_\nu}$. By \cite[Lemma 14]{A-preNichols} $\pi_i$ annihilates every term of the PBW-basis except $x_{\beta'_{M+1}}^{m}=x_i^m$, $m\ge 0$. Then
\begin{align*}
T_i^{\bq}\left( \iota_i\mathcal S^{-1}\pi_i \big( x_{\beta'_{M+1}}^{h_{M+1}}  x_{\beta'_{M}}^{h_{M}} \dots  x_{\beta'_{2}}^{h_{2}} \big)
\right) = 
\delta_{h_M,0}\dots \delta_{h_2,0}(-1)^{h_{M+1}} \,  T_i^{\bq}(x_i^{h_{M+1}}).
\end{align*}
Also, if $h_{M+1}>0$, then $T_i^{\bq}(x_i^{h_{M+1}}) \in \cU_\bq^{\le 0}$ by \eqref{eq:Ti-root-vectors}. 
Thus $x_{\mu}^{N^{\bq}_\mu} \otimes x_{\nu}^{N^{\bq}_\nu}$ can appear in \eqref{eq:coprod-x_gamma} only for those summands of \eqref{eq:xgammaprima-N-PBW} such that the middle part is $1$: that is, we only consider the terms  of \eqref{eq:xgammaprima-N-PBW} of the form
\begin{gather*}
x_{\mu'}^{N^{\bp}_{\mu'}} \ot 1 \ot x_{\beta'_{M+1}}^{h_{M+1}}  x_{\beta'_{M}}^{h_{M}} \dots x_{\beta'_{2}}^{h_{2}}.
\end{gather*}

Finally the right hand side of \eqref{eq:xgammaprima-N-PBW} must be $x_{\nu'}^{N^{\bq}_{\nu'}}$ by \eqref{eq:Ti-root-vectors}, so $c$ is the coefficient of $x_{\mu'}^{N^{\bq}_{\mu'}} \otimes x_{\nu'}^{N^{\bq}_{\nu'}}$ in $\varDelta(x_{\gamma'}^{N^{\bq}_{\gamma'}})$, which is non-zero by hypothesis.
\qed

\medbreak
We come back to the starting point of $\wmu$, $\wnu$, $\wgamma \in \wfO^\bq_+$.
By Theorem \ref{thm:cuntz-heck}
there exist $\bp\in\cX$,  $w\in\Hom (\bq,\bp)$ such that
\begin{align*}
\mu''&=w^{-1}(\mu), \, \nu'':= w^{-1}(\nu)\in\Delta_+^{\bp} \in \Z_{\geq 0} \alpha_1  + \Z_{\geq 0}\alpha_2;
\end{align*}
thus also $\gamma''=w^{-1}(\gamma)  \in \Z_{\geq 0} \alpha_1  + \Z_{\geq 0}\alpha_2$. 

\begin{claim}\label{claim:tres}
$[\xi_{\wmu''},\xi_{\wnu''}] \neq 0$.
\end{claim}

\noindent \emph{Proof of the Claim.} We follow the notation of \S \ref{subsubsec:subalgebras-distinguished}. Let $\rgo=(p_{ij})_{i,j\in\I_2}$ be the  $2\times2$-submatrix  
of $\bp=(p_{ij})_{i,j\in\I}$ corresponding to $\Jb=\I_2$. By Proposition \ref{prop:subalgebra-dist-pre-Nichols} there exists a  $\Z^{\I}$-graded Hopf algebra map $\widetilde{\varPhi}:Z_{\rgo}^+\to Z_{\bp}^+$. Hence $\widetilde{\varPhi}^t:\az_{\bp}\to \az_{\rgo}$ is a $\Z^{\I}$-graded Hopf algebra map which restricts to a $\Z^{\I}$-graded Lie algebra map $\widetilde{\varPhi}^t:\n^{\bp}\to \n^{\rgo}$. Let 
\begin{align*}
\n^{\bq}_{\Jb}:= 
\bigoplus_{\beta\in \Z_{\geq 0} \alpha_1  + \Z_{\geq 0}\alpha_2} (\n^\bq)_{\beta}.
\end{align*}
By Proposition \ref{prop:subalgebra-dist-pre-Nichols} \ref{item:subalgebra-dist-pre-Nichols-iii} 
$(\widetilde{\varPhi}^t)_{|\n^{\bq}_{\Jb}}:\n^{\bq}_{\Jb}\to \n^{\rgo}$ is injective. By Remark \ref{rem:nq-graded}, the non-trivial homogeneous components of $\n^{\rgo}$ are one-dimensional, so
\begin{align*}
(\n^{\rgo})_{\wmu''} &=\ku \, \widetilde{\varPhi}^t(\xi_{\wmu''}), & (\n^{\rgo})_{\wnu''} &=\ku \, \widetilde{\varPhi}^t(\xi_{\wnu''}).
\end{align*}

By \cite{AAR2}, $\n^{\rgo}$ is the positive part of a rank two semisimple Lie algebra. Hence \cite[\S 8.4]{Hum-libro}
applies:
\begin{align*}
[(\n^{\rgo})_{\wmu''},(\n^{\rgo})_{\wnu''}] = (\n^{\rgo})_{\gamma''}.
\end{align*}
Then $\widetilde{\varPhi}^t \big( [\xi_{\wmu''},\xi_{\wnu''}] \big) = \big[ \widetilde{\varPhi}^t(\xi_{\wmu''}), \widetilde{\varPhi}^t(\xi_{\wnu''}) \big] \neq 0$, and Claim \ref{claim:tres} follows. \qed

\medbreak

Pick a reduced expression $w =\sigma_{i_1}^{\bq}\dots \sigma_{i_t}$. 
By \cite[Lemma 8 (i)]{HY} there exists a reduced expression 
$w_0=\sigma_{i_1}^{\bq}\dots \sigma_{i_M}$ of $w_0$ which extends the reduced 
expression of $w$. By \eqref{eq:betak} there exist $j,k,\ell \in \I_M$ such that
\begin{align*}
\mu &= \sigma_{i_1}^{\bq}\dots \sigma_{i_{j-1}}(\alpha_{i_j}), &
\gamma &= \sigma_{i_1}^{\bq}\dots \sigma_{i_{k-1}}(\alpha_{i_k}), &
\nu &= \sigma_{i_1}^{\bq}\dots \sigma_{i_{\ell-1}}(\alpha_{i_\ell}).
\end{align*}
Up to exchange $\mu$ and $\nu$ we may assume that $j<l$, and necessarily $j<k<\ell$ because $N^{\bq}_{\mu}\mu +N^{\bq}_{\nu}\nu= N^{\bq}_{\gamma}\gamma$.

By Claim \ref{claim:tres} $[\xi_{\wmu''},\xi_{\wnu''}] \neq 0$ in $\n^{\bp}$.
As $\mu''=w^{-1}(\mu)$ is a positive root, we have that $t\le j$, see \cite[Lemma 8 (iii)]{HY}. 
Hence we may apply Claim \ref{claim:uno} repeatedly, starting with $\mu'',\nu'',\gamma''\in\wfO^{\bp}_+$, to conclude that  \eqref{eq:alpha-beta-gamma} holds.
\epf

In conclusion, Theorem \ref{th:determination-nq} follows by Lemma \ref{lem:graded-Lie-algebra-positive-roots} since the hypotheses \ref{item:graded-Lie-algebra-positive-roots-i}
and \ref{item:graded-Lie-algebra-positive-roots-ii} are satisfied by Remark \ref{rem:nq-graded} and Lemma \ref{lemma:generators-nq} respectively.

\section{Determination of $\n^\bq$}\label{sec:determination}
Finally  we explain how to obtain the Tables \ref{table:algebra-cartan}, \ref{table:algebra-super}, \ref{table:algebra-modular}, \ref{table:algebra-super-modular} and \ref{table:algebra-ufo}.

\subsection{The strategy}
By Theorem \ref{th:determination-nq}, we just need to compute the Cartan matrix $(a_{ij})$ of $\g^\bq$; clearly $a_{ii} = 2$. 
We proceed as follows:

\begin{enumerate}[leftmargin=*]
\item If $\wfO^{\bq}=\emptyset$, then set $\g^{\bq} := 0$ (this happens only once). \newline So we assume that $\vert \wfO^{\bq} \vert > 0$.

\item We compute the set $\varPi^{\bq}$ in \eqref{eq:Oq-simple} and fix a numeration $(\varpi_i)_{i \in \I_n}$.

\item Given $ \varpi_i \neq \varpi_j \in \varPi^{\bq}$, the $(i,j)$ entry is determined by
\begin{align}\label{eq:serre}
a_{ij} := - \sup \{ m\in \N_0: m \varpi_{i} + \varpi_{j} \in \wfO^{\bq}_+ \} \in \Z_{\leq 0}.
\end{align}

\end{enumerate}

We then go through all the diagrams of the classification in \cite{H-classif RS}
as described in \cite{AA17}. 
This is  now straightforward and we only give some examples 
to illustrate how to implement the methodology.  
By Lemma \ref{lemma:root-system-weyl-equivalence} it is enough to deal with just one element in each Weyl-equivalence class.

Following \cite{AA17}, we use a simplified notation for the expression of the positive roots as linear combinations of simple roots. For example, we write $34^25^36^2$ instead of $\alpha_3+2\alpha_4+ 3\alpha_5+ 2\alpha_6$. 

\begin{example}
Let $\toba_{\bq}$ be a Nichols algebra of Cartan type with Cartan matrix $C=(c_{ij})_{i,j\in\I}$. 
The matrix $\bq$ is written in terms of a root of unity $q\in\G_N'$ \cite[\S 4]{AA17}, 
where $N\ge 2$ if the Dynkin diagram of $C$ is simply laced, $N\neq 2,4$ if it has a double arrow and $N\ge 4$ if $C$ is of type $G_2$.

\begin{enumerate}[leftmargin=*]
\item If $N$ is coprime with all $c_{ij}$, $i\ne j$, then $\wfO^{\bq}$ is the root system associated with $C$. In particular, when the Dynkin diagram of $C$ is simply laced.

\item If $N$ is not coprime with all $c_{ij}$, $i\ne j$, then $\wfO^{\bq}$ is the root system of the Langlands dual group of $C$.
\end{enumerate}

\medspace

Indeed, $\fO^{\bq}=\Delta^{\bq}$ is the root system of $C$: all the roots are of Cartan type. If the Dynkin diagram of $C$ is simply laced, then $q_{ii}=q$ for all $i\in\I$. By Lemma \ref{lem:cartan-roots-invariant} $q_{\beta}=q$ for all $\beta\in\Delta_+^{\bq}$, so
\begin{align*}
\wfO^{\bq} &= \left\{ N\beta \,: \, \beta\in\Delta^{\bq} \right\}.
\end{align*}

Now assume that the Dynkin diagram of $C$ is not simply laced.
We set $c=2$ if $C$ is of types $B_{\theta}$, $C_{\theta}$, $F_4$, and $c=3$ if $C$ is of type $G_2$. Then $q_{ii}=q$ if $\alpha_i$ is  
a short root, and $q_{ii}=q^c$ if $\alpha_i$ is a long root. By Lemma \ref{lem:cartan-roots-invariant} we have
\begin{itemize}[leftmargin=*]
\item $q_{\beta}=q$ for all $\beta\in\Delta_{+,s}^{\bq}$, the set of short positive roots,
\item $q_{\beta}=q^c$ for all $\beta\in\Delta_{+,l}^{\bq}$, the set of long positive roots.
\end{itemize}
If $c$ does not divide $N$, then $N_{\beta}=N$ for all $\beta\in\Delta_+^{\bq}$, so $\wfO^{\bq} = \left\{ N\beta \,: \, \beta\in\Delta^{\bq} \right\}$ again. Otherwise, $N=cM$ for some $M\ge 2$, and
\begin{align*}
\wfO^{\bq} &= \left\{ \pm N\beta \,: \, \beta\in\Delta^{\bq}_{+,s} \right\} 
\cup \left\{ \pm M\beta \,: \, \beta\in\Delta^{\bq}_{+,l} \right\}.
\end{align*}
Then $\wfO^{\bq}$ is the Langlands dual since we exchange long and short roots. This result coincides with \cite{Lentner}.
\end{example}

\begin{example}
Let $\toba_{\bq}$ be a Nichols algebra of super type.
That is, the Weyl groupoid of $\bq$ is the one of a finite-dimensional contragredient Lie superalgebra $\g$ over a field of characteristic zero \cite{K-super}.
The matrix $\bq$ is written in terms of a root of unity $q\in\G_N'$, $N\ge 3$ \cite[\S 5]{AA17}.

\begin{enumerate}[leftmargin=*]
\item If $\bq$ is of type either $\supera{k-1}{\theta-k}$, $\superda{\alpha}$ or else $\superg$, then $\wfO^{\bq}$ is the root system of the Lie algebra $\g_0$ (the even part of $\g$). The same happens if $N$ is odd and $\bq$ is of type $\superb{k}{\theta-k}$, $\superd{k}{\theta-k}$ or $\superf$.

\item Assume that $N$ is even and $\bq$ is of type $\superb{k}{\theta-k}$, $\superd{k}{\theta-k}$ or $\superf$.
The Lie algebra $\g_0$ is a product of two simple Lie algebras and $\wfO^{\bq}$ is the root system of the product of one of these simple Lie algebra with the Langlands dual of the other, see Table \ref{table:algebra-super} for the precise description.
\end{enumerate}

\medbreak

For example, let $\toba_{\bq}$ be a Nichols algebra of super type $\superb{k}{\theta-k}$, $k\in\I$ \cite[\S 5.2]{AA17}.
Because of Lemma \ref{lemma:root-system-weyl-equivalence} we may choose just one of the matrices $\bq$ in \emph{loc. cit.}
Concretely we consider the diagram
\begin{align}\label{eq:dynkin-Btheta-super}
&\xymatrix{ \overset{q^{-2}}{\underset{\ }{\circ}}\ar  @{-}[r]^{q^2}  &
\overset{q^{-2}}{\underset{\ }{\circ}}\ar@{.}[r] &
\underset{k}{\overset{{-1}}{\underset{\
}{\circ}}}
\ar@{.}[r] & \overset{\quad q^{2}}{\underset{\ }{\circ}} \ar  @{-}[r]^{q^{-2}}  &
\overset{\quad q}{\underset{\ }{\circ}}},
 &
&\text{where }q \in \G'_N, \, N\neq 4.
\end{align}

For $1\leq i< j\leq\theta$ we set $\alpha_{ii}=\alpha_i$, $\alpha_{ij}= \alpha_i+ \alpha_{i+1}+ \dots+ \alpha_j$ and $\beta_{ij}=\alpha_{i\theta}+\alpha_{j\theta}$.
Hence the set of positive Cartan roots is
\begin{align*}
\fO^{\bq}_+= & \{ \alpha_{ij}: 1\leq i\leq j<k \}\cup \{ \alpha_{ij}: k< i\leq j\leq \theta \}\cup \{ \alpha_{i\theta}: 1\leq i\leq k \}\\
&
\cup \{ \beta_{ij}: 1\leq i< j\leq k \} \cup \{ \beta_{ij}: k< i< j\leq \theta \}.
\end{align*}

\smallbreak

First we assume that $N$ is odd. In this case, 
\begin{align*}
N^{\bq}_{\beta}= \begin{cases}
2N & \beta = \alpha_{i\theta}, \, i \in \I_k, \\
N, & \text{otherwise}.
\end{cases}
\end{align*}
We enumerate $\varPi^{\bq}$ as follows:
\begin{align*}
\varpi_j &=N \, \alpha_j, \quad  j\in\I-\{k\}; &
\varpi_{k} &= 2N \, \alpha_{k\theta}.
\end{align*}
By direct computation, 
\begin{align*}
a_{ij}= \begin{cases}
-2, & (i,j)=(k-1,k), (\theta,\theta-1); \\
-1, & |i-j|=1, \, (i,j) \neq (k\pm 1,k), (k,k+1),(\theta,\theta-1); \\ 
0, & \text{either }|i-j|>1 \text{ or }\{i,j\}=\{k,k+1\}.
\end{cases}. 
\end{align*}
Hence $\wfO^{\bq}$ is of type 
$C_{k}\times  B_{\theta-k}$.

\smallbreak

Next we assume that $N$ is even, $N=2M$. In this case, 
\begin{align*}
N^{\bq}_{\beta}= \begin{cases}
N & \beta = \alpha_{i\theta}, \, i\in \I, \\
M, & \text{otherwise}.
\end{cases}
\end{align*}
We enumerate $\varPi^{\bq}$ as follows:
\begin{align*}
\varpi_j &=M \, \alpha_j, \quad  j\in\I-\{k,\theta\}; &
\varpi_{k} &= N \, \alpha_{k\theta}; &
\varpi_{\theta} &= N\alpha_{\theta}.
\end{align*}
By direct computation, 
\begin{align*}
a_{ij}= \begin{cases}
-2, & (i,j)=(k-1,k), (\theta-1,\theta); \\
-1, & |i-j|=1, \, (i,j) \neq (k\pm 1,k), (k,k+1),(\theta-1,\theta); \\ 
0, & \text{either }|i-j|>1 \text{ or }\{i,j\}=\{k,k+1\}.
\end{cases}. 
\end{align*}
Hence $\wfO^{\bq}$ is of type 
$C_{k}\times  C_{\theta-k}$.
\end{example}

\begin{example}
Let $\toba_{\bq}$ be a Nichols algebra of type $\g(2, 6)$ \cite[\S 8.7]{AA17}.
Because of Lemma \ref{lemma:root-system-weyl-equivalence} we may choose just one of the matrices $\bq$ in \emph{loc. cit.}
Concretely we consider the diagram
\begin{align}\label{eq:dynkin-g(2,6)}
& \xymatrix{ \overset{\zeta}{\underset{\ }{\circ}}\ar  @{-}[r]^{\ztu}  &
\overset{\zeta}{\underset{\ }{\circ}} \ar  @{-}[r]^{\ztu}  & \overset{-1}{\underset{\ }{\circ}}
\ar  @{-}[r]^{\ztu}  & \overset{\zeta}{\underset{\ }{\circ}} \ar  @{-}[r]^{\ztu}  &
\overset{\zeta}{\underset{\ }{\circ}}} &
&\text{where }\zeta \in \G'_3.
\end{align}

Here, the set of positive Cartan roots is
\begin{align*}
\fO^{\bq}_+= & \{1, 12, 2,  12^23^245,  12^23^24,  123^24, 23^24, 4,12^23^24^25,
\\ 
& \quad  123^24^25,   23^24^25, 12^23^245, 123^245, 23^245, 45, 5\}.
\end{align*}
In this case $N^{\bq}_\beta=3$ for all $\beta\in\fO^{\bq}$. Now
\begin{align*}
\varPi^{\bq} = \{ \varpi_1 =2^3, \varpi_2 =1^3, \varpi_3 =2^33^64^3, \varpi_4 =5^3, \varpi_5 = 4^3 \}
\end{align*}
By direct computation, $a_{ij}= 
\begin{cases}
-1, & |i-j|=1, \\ 0, & |i-j|>1.
\end{cases}$. Hence $\wfO^{\bq}$ is of type $A_5$.
\end{example}

\begin{example}
Let $\toba_{\bq}$ be a Nichols algebra of unidentified type $\ufo(2)$ \cite[\S 10.2]{AA17}.
Again by Lemma \ref{lemma:root-system-weyl-equivalence} we   just consider the matrix $\bq$ with diagram
\begin{align}\label{eq:dynkin-ufo(5)}
& \xymatrix@C-6pt{\overset{\zeta}{\underset{\ }{\circ}}\ar
@{-}[r]^{\overline{\zeta}}  & \overset{\zeta}{\underset{\ }{\circ}}\ar  @{-}[r]^{\overline{\zeta}}  &
\overset{\zeta}{\underset{\ }{\circ}}
\ar  @{-}[r]^{\overline{\zeta}}  & \overset{\zeta}{\underset{\ }{\circ}}\ar  @{-}[r]^{\overline{\zeta}}  &
\overset{-1}{\underset{\ }{\circ}}\ar  @{-}[r]^{\overline{\zeta}}  &
\overset{\zeta}{\underset{\ }{\circ}}} &
&\text{where }\zeta \in \G'_4.
\end{align}

Here, the set of positive Cartan roots is
\begin{align*}
\fO^{\bq}_+= & \{1, 12, 2, 123, 23, 3, 1234, 234,34, 4, 12^23^34^35^36, 12^23^24^35^36, \\
&12^23^24^25^36, 123^24^25^36, 23^24^35^36, 23^24^25^36, 1234^25^36, 234^25^36, \\
& 34^25^36, 1^22^33^44^55^66^3, 12^33^44^55^66^3, 12^23^44^55^66^3, 12^23^34^35^36^2,\\
&  12^23^34^55^66^3, 12^23^24^35^36^2, 12^23^34^45^66^3, 12^23^24^25^36^2,123^24^35^36^2, \\
& 123^24^25^36^2, 1234^25^36^2, 23^24^35^36^2, 23^24^25^36^2,
234^25^36^2, 34^25^36^2, 6\}.
\end{align*}
In this case $N^{\bq}_\beta=4$ for all $\beta\in\fO^{\bq}$. Now
\begin{align*}
\varPi^{\bq} = \{ \varpi_1 =1^4, \varpi_2 =2^4, \varpi_3 =3^4, \varpi_4 =4^4, \varpi_5 = 3^44^85^{12}6^4, \varpi_6 =6^4 \}
\end{align*}
By direct computation, 
\begin{align*}
-1 &= a_{12} =a_{23} =a_{34} =a_{25} =a_{56} =a_{21} =a_{32} =a_{43} =a_{52} =a_{65},
\end{align*}
and $a_{ij}=0$ for the remaining pairs $i\neq j\in\I_6$. Hence $\wfO^{\bq}$ is of type $E_6$.
\end{example}

\appendix
\section{Presentations for type $B_2$}
\subsection{Overview}

One of the Referees kindly asked us to illustrate the paper by  presenting explicitly the different algebras appearing in \eqref{eq:Lusztig-exact-sequence},
\eqref{eq:DCKP-exact-sequence}, \eqref{eq:intro-pre-Nichols-exact-sequence} and \eqref{eq:extension-braided-lu-intro}
when $\mathfrak g$ is of type $B_2$.
For completeness we report on diagrams of type $B_2$ that is one of the  following:

\begin{enumerate}[leftmargin=*,label=\rm{(\alph*)}]
\item Cartan $B_2$, $q \in \G'_N$,  $N$ odd; diagram $\Dchaintwo{q}{q^{-2}}{q^2}$. See \cite[4.2.5]{AA17} for details.

\item Cartan $B_2$, $q \in \G'_N$,  $N = 2M$ even; same diagram, see \cite[4.2.5]{AA17}.

\item Standard $B_2$, $\zeta \in \G_3'$, $\mathbb{J}=\emptyset$, $\Dchaintwo{\zeta }{-\zeta}{-\zeta^{-1}}$.
See \cite[6.1.4]{AA17} for details.

\item  Standard $B_2$, $\zeta \in \G_3'$, $\mathbb{J}=\{1\}$, $\Dchaintwo{\zeta }{-\zeta}{-1}$.
See \cite[6.1.4]{AA17} for details.
\end{enumerate}

\subsection{The algebras $\mathcal{U}_q(\mathfrak g)$ and $\mathfrak{u}_q(\mathfrak g)$ in \eqref{eq:Lusztig-exact-sequence}}
Let $\g$ be a simple, finite-dimensional, complex Lie algebra.
The Hopf algebra $\cU_q(\g)$ is the \emph{quantum group at a root of 1} 
introduced in \cite{L - mod rep} when $\g$ is simply laced and the order of $q$ is odd.
The presentation by generators and relations is given in \cite[Subsection 2.3]{L-fdHa-JAMS}.
The definition and presentation were extended to $\g$ multiply laced in \cite{L-qgps-at-roots}, still with the order of $q$ odd
and prime to $3$ if $\g$ is of type $G_2$. 

\medbreak Since the well-known algebras $\mathcal{U}_q(\mathfrak g)$ and $\mathfrak{u}_q(\mathfrak g)$
appeared long ago,
we comment on the general shape of the relations and give precise references.
The algebra $\cU_q(\g)$ has a triangular decomposition 
\begin{align*}
\cU_q(\g) \simeq \cU_q^-(\g) \otimes \cU_q^0(\g) \otimes \cU_q^+(\g)
\end{align*}
where $\cU_q^+(\g)$ is isomorphic to what we call the Lusztig algebra $\luq{\bq}$ for a suitable $\bq$
and analogously for $\cU_q^-(\g)$. See   the presentation below.
In turn $\cU_q^0(\g) \simeq \ku \Gamma \otimes U(\mathfrak h)$ where $\Gamma$ is a finite group 
and $\mathfrak h$ is an abelian Lie algebra. Besides the relations from the subalgebras $\cU_q^-(\g)$, $\cU_q^0(\g)$ and $\cU_q^+(\g)$,
there are relations accounting for the actions of $\cU_q^0(\g)$ on $\cU_q^{\pm}(\g)$ and the relations of the kind $[E_i, F_j] = \delta_{ij}\cdots$.

\medbreak
The Hopf algebra $\ug_q(\g)$ is the so called \emph{Frobenius-Lusztig kernel} or also \emph{small quantum group}.
It was introduced in \cite[Section 5]{L-fdHa-JAMS} when $\g$ is simply laced and the order of $q$ is odd
and in \cite{L-qgps-at-roots} when $\g$ is multiply laced and the order of $q$ is odd and prime to $3$ if $\g$ is of type $G_2$.

The algebra $\ug_q(\g)$ has a triangular decomposition 
\begin{align*}
\ug_q(\g) \simeq \ug_q^-(\g) \otimes \ug_q^0(\g) \otimes \ug_q^+(\g)
\end{align*}
where $\ug_q^+(\g)$ is isomorphic to the Nichols algebra $\toba_{\bq}$ for a suitable $\bq$
and analogously for $\ug_q^-(\g)$. In turn $\ug_q^0(\g) \simeq \ku \Gamma$ where $\Gamma$ is the same finite group as above.
The relations of $\ug_q(\g)$  arise from the subalgebras $\ug_q^-(\g)$, $\ug_q^0(\g)$ and $\ug_q^+(\g)$,
from the actions of $\ug_q^0(\g)$ on $\ug_q^{\pm}(\g)$ and again those of the kind $[E_i, F_j] = \delta_{ij}\cdots$.
Variations of these Hopf algebras for $q$ of even order appear in \cite{Lentner,L-libro,Negron,Sawin}.

\subsection{The algebras  $U_q(\mathfrak g)$ and $\mathcal{O}(G^d)$ in  \eqref{eq:DCKP-exact-sequence}}
The Hopf algebra $U_q(\g)$ is the \emph{quantum group at a root of 1} 
introduced by generators and relations in \cite{DK} when the order of $q$ is odd and prime to $3$ if $\g$ is of type $G_2$.
We do not know where the variation for even order appeared the first time.
Since the algebra $U_q(\mathfrak g)$ appeared long ago,
we just comment on the general shape of the relations.
The algebra $U_q(\g)$ has a triangular decomposition 
\begin{align*}
U_q(\g) \simeq U_q^-(\g) \otimes U_q^0(\g) \otimes U_q^+(\g)
\end{align*}
where the relations of $U_q^+(\g)$ are just the quantum Serre relations
and analogously for $U_q^-(\g)$. 
In turn $U_q^0(\g) \simeq \ku \Gamma$ where 
$\Gamma$ is now a  free abelian group of finite rank (= the rank of $\g$).
Besides the quantum Serre relations from the subalgebras $U_q^{\pm}(\g)$, 
there are relations accounting for the actions of $U_q^0(\g)$ on $U_q^{\pm}(\g)$ and the relations of the kind $[E_i, F_j] = \delta_{ij}\cdots$.

\medbreak
The  central Hopf subalgebra $\cO (G^d)$ appeared in \cite{DKP} and was studied in detail in \cite{DP}. 
As an algebra is a polynomial algebra in $\dim \g$-variables. The algebraic group $G^d$ is solvable.

\bigbreak
\begin{table}[ht]
	\setlength{\unitlength}{1mm}
	\begin{tabular}{l|l|l}
		& Cartan $B_2$, $N$ odd &Cartan $B_2$, $N=2M$ even 
		\\
		\hline \hline
		\rule{0mm}{16pt}
		$\toba_{\bq}$ & $\ku \langle x_1, x_2| \, x_{122},x_{1112},x_\alpha^N\rangle$ & $\ku \langle x_1, x_2| \, x_{122},x_{1112},x_2^M,x_{112}^M,x_1^N,x_{12}^N\rangle$ \\
		\hline
		\rule{0mm}{16pt}
		$\dpn_{\bq}$ & $\ku \langle x_1, x_2| \, x_{122},x_{1112}\rangle$& $\ku \langle x_1, x_2| \, x_{122},x_{1112}\rangle$\\
		\hline
		\rule{0mm}{16pt}
		$Z_{\bq}^+$ & $\ku \langle x_1^N, x_2^N, x_{12}^N, x_{112}^N | \, x_\alpha^Nx_\beta^N-x_\beta^N x_\alpha^N \rangle$&
		$\ku \langle x_1^N, x_2^M, x_{12}^N, x_{112}^M | \, x_\alpha^{N_\alpha}x_\beta^{N_\beta} -x_\beta^{N_\beta} x_\alpha^{N_\alpha} \rangle$ \\
		\hline
		\rule{0mm}{16pt}
		$\az_{\bq}$ & $\ku \langle \xi_1, \xi_{12}, \xi_{112},\xi_2 | \, [\xi_1,\xi_{112}] - [\xi_2,\xi_{12}],$ 
		& $\ku \langle \xi_1, \xi_{12}, \xi_{112},\xi_2 | \, [\xi_1,\xi_{112}] - [\xi_2,\xi_{12}],$ \\
		& $\quad [\xi_1,\xi_2]-(1 - \tfrac{1}{q^2})^N \xi_{12},$ 
		&  $ \quad [\xi_2,\xi_1]-(1- \tfrac{1}{q})^M (1 - \tfrac{1}{q^2})^M q_{21}^{M(M-1)}\xi_{112}$, \\
		& $\quad [\xi_{12},\xi_{1}]-(1- \tfrac{1}{q})^N(1 - \tfrac{1}{q^2})^N \xi_{112} \rangle$
		&$\quad [\xi_{112},\xi_2]- (1 - \tfrac{1}{q^2})^M q_{21}^{M^2}\xi_{12} \rangle$
		\\
		\hline \hline
		\rule{0mm}{16pt}
		& Standard $B_2$, $\zeta \in \G_3'$, $\mathbb{J}=\{1\}$ & Standard $B_2$, $\zeta \in \G_3'$, $\mathbb{J}=\emptyset$
		\\
		\hline \hline
		\rule{0mm}{16pt}
		$\toba_{\bq}$& $\ku\langle x_1, x_2 | \, [x_{112}, x_{12}]_c, x_1^3, x_2^2 \rangle$ & $\ku\langle x_1, x_2| \, x_{122}, x_1^3, x_2^6, x_{112}^6 \rangle$\\
		\hline
		\rule{0mm}{16pt}
		$\dpn_{\bq}$ &$\simeq \toba_\bq\simeq \lu_\bq$ & $\ku\langle x_1, x_2| \, x_{122}, x_1^3 \rangle$ \\
		\hline
		\rule{0mm}{16pt}
		$Z_{\bq}^+$ & $\ku$ &$\ku\langle x_2^6, x_{112}^6| \, x_2^6x_{112}^6-x_{112}^6x_2^6\rangle$ \\
		\hline
		\rule{0mm}{16pt}
		$\az_{\bq}$ & $\ku$ & $\ku\langle \xi_2,\xi_{112}| \, [\xi_2,\xi_{112}] \rangle$ \\
		
	\end{tabular}
	\\[2mm] \
	\caption{Generators and relations of some algebras of type $B_2$.}
	\label{table:cartan}
\end{table}

\subsection{The algebras $\toba_{\bq}$, $\dpn_{\bq}$ and $Z_{\bq}^+$ in \eqref{eq:intro-pre-Nichols-exact-sequence} and $\az_{\bq}$ in  \eqref{eq:extension-braided-lu-intro}}
We give the generators and relations of these algebras in Table \ref{table:cartan}.

\medbreak
Assume first that we are in  Cartan type with $N$ odd. 
As we said the algebras $\toba_{\bq}$ were introduced in \cite{L-fdHa-JAMS,L-qgps-at-roots}
as the positive parts of the small quantum groups; 
where the defining relations  (quantum Serre relations and powers of root vectors)
appeared first is difficult to track back but they were well-known by the end of the 90's, see e.~g. \cite[Theorem 4.6]{AS-pointed} and references therein. 
For type $B_2$, the presentation is explicitly stated in \cite{BDR} .

\medbreak
We also said that the relations of $\dpn_{\bq}$ in this case are simply the quantum Serre relations, hence
$Z_{\bq}^+$ is a polynomial algebra spanned by the powers of the root vectors. The underlying algebraic group was determined in \cite{DP}
and is related to the group $G^d$. 
The algebra $\az_{\bq}$ appeared in \cite{L-fdHa-JAMS,L-qgps-at-roots} as the enveloping algebra of the positive part of $\g$.

\medbreak
The algebras $\toba_{\bq}$ for  Cartan type with $N$ even are described in \cite{L-libro}.
We do not know the first appearence of the relations; 
certainly they are present in \cite{Ang-crelle,A-preNichols,AAR2} for $\toba_{\bq}$, $\dpn_{\bq}$ and $Z_{\bq}^+$, and $\az_{\bq}$ respectively.

\medbreak
For  standard type, the existence of $\toba_{\bq}$ follows by the general theory of Nichols algebras and the finite-dimensionality follows from \cite{H-classif RS}. The relations were described first in \cite{Ang-ant}. 
The algebras $\dpn_{\bq}$ and $Z_{\bq}^+$ appeared for the first time, together with their presentations,
in \cite{A-preNichols}. The algebra $\az_{\bq}$ was introduced in \cite{AAR2}.

\subsection{The algebras $\lu_{\bq}$ in  (1.5)}
We give below the generators and relations of the algebras $\lu_{\bq}$.
As we said, the presentation was given in \cite{L-qgps-at-roots} for Cartan type, $N$ odd. 
For standard type, and probably for Cartan type with $N$ even, the presentation appears first in \cite{AAR1}.

\medbreak
Let $\Delta^+ = \{1,2,12,112\}$.

\subsubsection{The algebra $\lu_{\bq}$, Cartan $B_2$, $q \in \G'_N$, $N$ odd}
The algebra 
$\lu_{\bq}$ is presented by generators $y_{\alpha}^{(n)}$, $\alpha \in \Delta^+$, $n\in \N$,
with defining relations
\begin{align}\label{eq:luq-Cartan-Nodd-1}
&y_\alpha^{(n)} y_\alpha^{(m)} - \binom{n+m}{n}_q y_\alpha^{(n+m)},  \qquad n,m\in \N, \quad \alpha \in \Delta^+,
\\ \label{eq:luq-Cartan-Nodd-2}
&[y_{2}^{(N)}, \  y_{112}^{(N)}]_c,  
\\ \label{eq:luq-Cartan-Nodd-3}
&[y_{2}^{(n)}, \  y_{12}^{(m)}]_c, \qquad \qquad \qquad \qquad n,m\in \N,
\\ \label{eq:luq-Cartan-Nodd-4}
&[y_{112}^{(n)}, \ y_1^{(m)}]_c, \qquad \qquad\qquad \qquad n,m\in \N,
\\\label{eq:luq-Cartan-Nodd-5}
&[y_{112}^{(n)}, \ y_{12}^{(m)}]_c,  \qquad \qquad\qquad \qquad n,m\in \N,
\\ \label{eq:luq-Cartan-Nodd-6}
&[y_2,y_1]_c - (1 - \tfrac{1}{q^2})y_{12}, 
\\ \label{eq:luq-Cartan-Nodd-7}
&[y_{12},y_1]_c - q(1 - \tfrac{1}{q^2})y_{112}, 
\\ \label{eq:luq-Cartan-Nodd-8}
&[y_{2},y_{112}]_c - q(1 - \tfrac{1}{q^2})(1+q)y^{(2)}_{12},
\\ \label{eq:luq-Cartan-Nodd-9}
&[y_{12}^{(N)},y_1]_c - q_{21}^{N-1}(1 - \tfrac{1}{q^2}) y_{112}y_{12}^{(N-1)},
\\\label{eq:luq-Cartan-Nodd-10}
&[y_{2},y_{112}^{(N)}]_c - (1+q^{-1})q_{21}^{2N-2}(1 - \tfrac{1}{q^2}) y_{112}^{(N-1)}y_{12}^{(2)} 
\\ \label{eq:luq-Cartan-Nodd-11}
&\begin{aligned}
&[y_2^{(N)},y_1^{(N)}]_c  - (1 - \tfrac{1}{q^2})y_{12}^{(N)} 
\\&- \sum_{k=1}^{N-1} (1- \tfrac{1}{q})^kq_{21}^{\frac{2N^2-k^2-k}{2}} y_1^{(N-k)}y_{12}^{(k)}y_2^{(N-k)}, 
\end{aligned}
\\ \label{eq:luq-Cartan-Nodd-12}
&\left\{\begin{aligned}
&[y_{12}^{(N)},y_1^{(N)}]_c -  q(1 - \tfrac{1}{q^2})y_{112}^{(N)} \\
&- \sum_{k=1}^{N-1} (1 - \tfrac{1}{q^2})^kq^{-k(k-1)}q_{21}^{\frac{2N^2-k^2-k}{2}} y_{1}^{(N-k)}y_{112}^{(k)}y_{12}^{(N-k)},
\end{aligned}\right.
\\\label{eq:luq-Cartan-Nodd-13}
&[y_2^{(N)},y_1]_c  -   q_{21}^{N-1}(1 - \tfrac{1}{q^2}) y_{12}y_2^{(N-1)}, 
\\\label{eq:luq-Cartan-Nodd-14}
&[y_2,y_1^{(N)}]_c - qq_{21}^{N-2}(1 - \tfrac{1}{q^2})y_1^{(N-2)} (-q_{21} y_1y_{12}+ (1- \tfrac{1}{q})y_{112}),
\\\label{eq:luq-Cartan-Nodd-15}
&[y_{12},y_1^{(N)}]_c - q_{21}^{N-1}(1 - \tfrac{1}{q^2}) y_1^{(N-1)}y_{112}, 
\\ \label{eq:luq-Cartan-Nodd-16}
&[y_{2}^{(N)},y_{112}]_c   - (1+q^{-1})q_{21}^{2N-2}(1 - \tfrac{1}{q^2}) y_{12}^{(2)}y_{2}^{(N-1)}.
\end{align}

\subsubsection{The algebra $\lu_{\bq}$, Cartan $B_2$, $q \in \G'_\N$, $N = 2M$ even}
The algebra 
$\lu_{\bq}$ is presented by generators $y_{\alpha}^{(n)}$, $\alpha \in \Delta^+$, $n\in \N$,
with defining relations \eqref{eq:luq-Cartan-Nodd-1}, \eqref{eq:luq-Cartan-Nodd-3}, \eqref{eq:luq-Cartan-Nodd-4}, \eqref{eq:luq-Cartan-Nodd-5}, 
\eqref{eq:luq-Cartan-Nodd-6}, \eqref{eq:luq-Cartan-Nodd-7},  \eqref{eq:luq-Cartan-Nodd-8}, \eqref{eq:luq-Cartan-Nodd-9}, 
\eqref{eq:luq-Cartan-Nodd-14}, \eqref{eq:luq-Cartan-Nodd-15}, 

\begin{align}\label{eq:luq-Cartan-Neven-1}
&[y_{12}^{(N)},\   y_1^{(N)}]_c, 
\\\label{eq:luq-Cartan-Neven-2}
&[y_{2}^{(M)},y_{112}]_c - (1+q^{-1})q_{21}^{2M-2}(1 - \tfrac{1}{q^2}) y_{12}^{(2)}y_{2}^{(M-1)},
\\ \label{eq:luq-Cartan-Neven-3}
&[y_{2},y_{112}^{(M)}]_c - (1+q^{-1})q_{21}^{2M-2}(1 - \tfrac{1}{q^2}) y_{112}^{(M-1)}y_{12}^{(2)}, 
\\\label{eq:luq-Cartan-Neven-4}
&[y_2^{(M)},y_1]_c - q_{21}^{M-1}(1 - \tfrac{1}{q^2}) y_{12}y_2^{(M-1)}, 
\\ \label{eq:luq-Cartan-Neven-5}
&\left\{\begin{aligned}
&[y_{2}^{(M)},y_{112}^{(M)}]_c - (1 - \tfrac{1}{q^2})^M q_{21}^{M^2} y_{12}^{(N)} 
\\&- \sum_{k=1}^{M-1} q^{4k^2} (2)_q^{2k} (1 - \tfrac{1}{q^2})^{2k} q_{21}^{4M+4k^2-2k} y_{112}^{(M-k)}y_{12}^{(2k)}y_2^{(M-k)} , 
\end{aligned}\right.
\\ \label{eq:luq-Cartan-Neven-6}
&\left\{\begin{aligned}
&[y_2^{(M)},y_1^{(N)}]_c - (1- \tfrac{1}{q})^M(1 - \tfrac{1}{q^2})^M q_{21}^{M^2-M} y_{112}^{(M)} 
\\&- \sum_{k=1}^{M-1} (1- \tfrac{1}{q})^k(1 - \tfrac{1}{q^2})^k q_{21}^{2M^2-k^2-k} y_1^{(N-2k)}y_{112}^{(k)}y_2^{(M-k)} 
\\ 
&- \sum_{k=1}^{M-1} (1- \tfrac{1}{q})^kq_{21}^{\frac{N^2-k^2-k}{2}} y_1^{(N-k)}y_{12}^{(k)}y_2^{(M-k)}.
\end{aligned}\right.
\end{align}

\subsubsection{The algebra $\lu_{\bq}$, Standard $B_2$, $\zeta \in \G_3'$, $\mathbb{J}=\emptyset$}
The algebra $\lu_{\bq}$ is presented by generators $y_1^{(n)}, y_2^{(n)}$, $n\in \N$, with defining relations
\eqref{eq:luq-Cartan-Nodd-1}, \eqref{eq:luq-Cartan-Nodd-3}, \eqref{eq:luq-Cartan-Nodd-4}, \eqref{eq:luq-Cartan-Nodd-5}, 
\begin{align}\label{eq:luq-standard-1}
&y_1^{3}, 
\\\label{eq:luq-standard-2}
& [y_2^{(6)},y_1]_c- (1+\zeta)q_{21}^5 y_{12}y_2^{(5)}, 
\\\label{eq:luq-standard-3}
&[y_{12},y_1]_c-y_{112},
\\\label{eq:luq-standard-4}
&[y_2,y_{112}]_c-(1+\zeta)q_{12}^{-1}y_{12}^{(2)},
\\ \label{eq:luq-standard-5}
& [y_2,y_1]_c-(1+\zeta)y_{12}, 
\\\label{eq:luq-standard-6}
&[y_2^{(6)},y_{112}^{(6)}]_c- (1+\zeta)q_{21}^{11} y_{12}^{(2)}y_2^{(5)},
\\ \label{eq:luq-standard-7}
&[y_2^{(6)},y_{112}^{(6)}]_c- (1+\zeta)q_{21}^{11} y_{112}^{(5)}y_{12}^{(2)},
\\ \label{eq:luq-standard-8}
& [y_2^{(6)},y_{112}^{(6)}]_c-\sum_{k=1}^5(1+\zeta)^{12-2k}(-\zeta)^{k^2+k}q_{21}^{-2k^2-120} y_{112}^{(k)}y_{12}^{(12-2k)}y_2^{(k)}.
\end{align}

\end{document}